\documentclass[12pt]{amsart}
\usepackage{amsmath}
\usepackage{tgpagella}
\usepackage{euler}
\usepackage[T1]{fontenc}
\usepackage{amsmath, amssymb,mathabx}
\usepackage[hidelinks]{hyperref}
\usepackage[english]{babel}
\usepackage{mathrsfs}
\usepackage{eucal}
\usepackage[all]{xy}
\usepackage{enumitem}
\usepackage{imakeidx}

\newtheorem{thm}{Theorem}[section]
\newtheorem*{thm*}{Theorem}
\newtheorem{lem}[thm]{Lemma}
\newtheorem{fact}[thm]{Fact}
\newtheorem{prop}[thm]{Proposition}
\newtheorem*{prop*}{Proposition}

\theoremstyle{definition}
\newtheorem{defn}[thm]{Definition}

\newtheorem{example}[thm]{Example}

\newtheorem{exercise}[thm]{Exercise}
\newtheorem{comment}[thm]{Comment}

\def\e{\epsilon}

\def\bb{\mathbb}

\def\bb{\mathbb}

\def\cal{\mathcal}

\def\M{\mathcal{M}}
\def\N{\mathcal{N}}
\newcommand{\cstar}{$\mathrm{C}^*$}
\def\bbR{\mathbf R}

\newcommand\cU{{\cal U}}
\newcommand\cV{{\cal V}}
\newcommand\cL{{\cal L}}

\newcommand\bbN{{\mathbf N}}
\newcommand\bbC{\mathbf C}

\newcommand{\dotminus}{ 
\buildrel\textstyle\ .\over{\hbox{ 
\vrule height3pt depth0pt width0pt}{\smash-} 
}\ }

\makeindex
\begin{document}
\title{An introduction to continuous model theory}

\author{Bradd Hart}
\address{Department of Mathematics and Statistics, McMaster University, 1280 Main St., Hamilton ON, Canada L8S 4K1}
\email{hartb@mcmaster.ca}
\urladdr{http://ms.mcmaster.ca/~bradd/}
\thanks{The author was supported by the NSERC.}

\maketitle

\begin{abstract}
    We present an introduction to modern continuous model theory with an emphasis on its interactions with topics covered in this volume such as \cstar-algebras and von Neumann algebras.  The role of ultraproducts is highlighted and expositions of definable sets, imaginaries, quantifier elimination and separable categoricity are included.
\end{abstract}

\section{Introduction}

It is not a new idea to have more truth values than just true and false. For example, \L ukasiewicz and Tarski suggested a real-valued logic in the 1930's \cite{LT}. The modern understanding of continuous model theory has other antecedents as well.  Chang and Keisler wrote a book \cite{CK} in which they consider a logic with truth values in arbitrary compact Hausdorff spaces. This logic had many quantifiers and did not catch on; see Ben Yaacov \cite{BY-quant} for further commentary.  Henson introduced positive bounded logic in \cite{Henson}, mostly in the service of the model theory of normed vector spaces with additional structure.  Ben Yaacov introduced the very general setting of cats (short for compact abstract theories) in \cite{BY-cats} in order to capture expansions of first order logic that would allow the seamless use of hyperimaginaries as introduced in \cite{HKP}; this line of research built on earlier work around Robinson theories by Hrushovski, Pillay and others. In particular, it was realized that Hausdorff cats interpreted a metric. Ben Yaacov, Berenstein, Henson and Usvyatsov then introduced a logic for metric structures which is the subject of this paper.  This logic is first presented in \cite{BYU}, which strangely was published after the de facto standard \cite{BBHU}.  
%The abstract model theory of metric structures is discussed in a paper of Iovino \cite{Iovino} in which he proves the Lindstr\"om-like characterization of continuous first order logic which is presented in the lectures.  Note that this paper pre-dates the definition of metric structure!

The history of the ultraproduct is very interesting and more colourful than most model theorists are aware.  A short history is contained in \cite{Sherman:ultra}.  The thumbnail sketch is that an ultraproduct-like construction in the operator algebra context was introduced by Kaplansky and Wright in the early 1950's.  Sakai also used an ultraproduct in the early 1960's, again in the service of operator algebras.  McDuff's systematic use of ultraproducts in her paper \cite{McDuff:1} highlighted how important the construction had become in the operator algebra world.  In parallel, the seminal introduction of the ultraproduct in model theory by \L o\' s \cite{Los} and its use by Robinson in the development of non-standard analysis were key early moments in model theory.  It does not seem that anyone at that time saw anything more than a formal connection between the two uses of ultraproducts.  Keisler \cite{Keisler1} provides a brief history on the model theory side in which he gives the nod for a precursor to the ultraproduct to Skolem and interestingly to Hewitt  \cite{Hewitt}, who was working on rings of continuous functions.  Although there was prior work done by Krivine, Stern and others, the first systematic connection between the analytic ultraproduct construction and model theory, particularly for Banach spaces, was the work of Henson \cite{Henson}, in which he also introduces his positive bounded logic.  A general exposition of this logic appears with Iovino in \cite{HenIov}.

The ultraproduct plays a slightly different role in the development of continuous model theory than it does in classical logic.  While the connection between the ultraproduct construction and the compactness theorem is the same for both logics, it is the exploratory nature of the ultraproduct which is different, by which we mean that there are classes of structures appearing in analysis where an ultraproduct is used but without a formal language being present.  The fact that mathematicians are using ultraproducts in certain settings "in the wild" often leads one to examine whether the class can in fact be presented in a model theoretic framework.  Examples relevant to this volume where this has been done successfully include Banach spaces, \cstar-algebras, tracial von Neumann algebras, and W$^*$-probability spaces.

As we will see in Example \ref{CFOL}, every classical elementary class can be seen as a continuous elementary class.  Continuous logic, however, is a \emph{proper} extension of classical first order logic.    While there are a number of ways to see this, here is one suggestion based on a Borel complexity argument due to Farah, Toms and T\"ornquist. In \cite{FTT}, they show that the isomorphism relation on the class of separable \cstar-algebras is turbulent and  not Borel-reducible to any class of countable classical structures (see their article for the relevant definitions).  Of course, the isomorphism relation for every class of countable classical structures is Borel-reducible to itself, so the isomorphism relation for the class of separable \cstar-algebras is not equivalent to any class of countable classical structures.  As we will see in this article, the class of \cstar-algebras is an elementary class (Example \ref{C*-axioms}) and by this complexity argument, this class is not equivalent to any classical elementary class.

We now give a brief outline of the paper.  After a quick reminder of the basics of ultrafilters in Section \ref{sec-prelim}, we introduce the syntax and semantics of continuous logic in Sections \ref{met-str} and \ref{sec-form}. 
 Basic continuous model theory is then developed in Section \ref{sec-mt}, where a generalized notion of formula is introduced. 
 As was indicated above, ultraproducts play a crucial role in continuous model theory and the ultraproduct of metric structures is introduced in Section \ref{sec-ultra}. Section \ref{sec-types} introduces the notion of type in continuous model theory and the type space is used to give another characterization of general formulas.  A critical section is Section \ref{sec-def}, where definable sets are discussed and contrasted with the notion of definable set from classical logic. Section \ref{sec-qe} highlights the formulation of quantifier elimination in continuous logic and gives several examples. In Section \ref{sec-imag}, imaginary elements in the continuous setting are considered and the conceptual completeness theorem is stated in this context.  We end the article with a discussion of omitting types in Section \ref{sec-omit} and separable categoricity in Section \ref{sec-sep-cat}. 

 Inevitably, there are topics that we left out.  In particular, we do not include a discussion of a proof system for continuous logic.  This topic is covered in \cite{BYP} and is important for the topic of decidability in continuous logic (see, for example, Section 3 of Goldbring's paper in this volume).  Another topic not discussed here is stability theory in continuous logic; one can find the basics of this theory in \cite{BBHU} and \cite{BYU}.

 I'd like to thank Isaac Goldbring for our many discussions regarding this article and for his patience as an editor.

\section{Preliminaries on ultrafilters and ultraproducts}\label{sec-prelim}

Ultrafilters and ultraproducts will play an important role in continuous model theory.  We cover the requisite facts in this section.
\begin{defn}
Suppose that $X$ is a set and $F \subset P(X)$, where $P(X)$ is the power set of $X$.
\begin{enumerate}
\item $F$ is a (proper) \emph{filter}\index{filter} on $X$ if:
\begin{enumerate}
\item $\emptyset \not \in F$, 
\item $A, B \in F$ implies $A \cap B \in F$, and 
\item $A \subset B \subset X$ and $A \in F$ implies $B \in F$.
\end{enumerate}
\item $F$ is an \emph{ultrafilter}\index{ultrafilter} on $X$ if it is a filter on $X$ and, for every $A \subset X$, either $A \in F$ or $X \setminus A \in F$ (but not both).
\end{enumerate}
\end{defn}

The following is an important characterization of ultrafilters.  Note that we are assuming the axiom of choice throughout this paper.
\begin{exercise}
Suppose $F\subset P(X)$.
\begin{enumerate}
    \item $F$ is an ultrafilter on $X$ if and only if $F$ is a maximal (under inclusion) filter on $X$.
    \item $F$ is contained in an ultrafilter on $X$ if and only if $F$ has the finite intersection property, that is, for any finitely many $A_1,\ldots, A_n \in F$, $\bigcap_{i = 1}^n A_i \neq \emptyset$.
\end{enumerate}
\end{exercise}

\begin{example}

\

\begin{enumerate}
\item If $X$ is any infinite set and $F$ is the collection of co-finite subsets of $X$, then $F$ is a filter on $X$, called the \emph{Fr\'echet filter}\index{Fr\'echet filter} on $X$. Note that the Fr\'echet filter is not an ultrafilter.
\item If $X$ is any set, $a \in X$ is any element, and $F$ is the collection of subsets of $X$ which contain $a$, then $F$ is an ultrafilter on $X$, called the \emph{principal ultrafilter on $X$ generated by $a$}\index{principal ultrafilter}.
\end{enumerate}
\end{example}

\begin{exercise}
If $X$ is an infinite set, then an ultrafilter on $X$ is non-principal if and only if it contains the Fr\'echet filter on $X$.  Consequently, non-principal ultrafilters on $X$ exist.
\end{exercise}

The notion of an ultralimit will be important throughout this paper. 

\begin{exercise}
Suppose that $(r_i : i \in I)$ is a sequence of real numbers and $\cU$ is an ultrafilter on $I$.  Further  suppose that there is a real number $B$ such that $\{ i \in I : |r_i| < B \} \in \cU$.  Prove that there is a unique real number $r$ satisfying, for every $\e >0$, that $\{ i \in I : |r_i - r| < \e \} \in \cU$.
\end{exercise}

\begin{defn}
In the context of the previous exercise, we say that $r$ is the \emph{ultralimit of $(r_i : i \in I)$\index{ultralimit} along $\cU$}, written $r = \lim_\cU r_i$.
\end{defn}

% \begin{exercise}
% Show that if the ultralimit exists, then it is unique.  Moreover, if \\$(r_i : i \in I)$ is a sequence of real numbers, $\cU$ is an ultrafilter on $I$, and there is a real number $B$ such that $\{ i \in I : |r_i| < B \} \in \cU$, show that the ultralimit of $(r_i : i \in I)$ along $\cU$ exists.
% \end{exercise}

We now introduce the notion of the metric ultraproduct for bounded metric spaces.  Suppose that $(X_i,d_i)$ is an $I$-indexed family of metric spaces that is \emph{uniformly bounded}\index{uniformly bounded}, that is, there is some $B>0$ such that the diameter of each $X_i$ is bounded by $B$ (meaning for all $i\in I$ and $a,b \in X_i$, $d_i(a,b) \leq B$). Fix an ultrafilter $\cU$ on $I$. Set $X = \prod_{i\in I} X_i$ and define, for $\bar a,\bar b \in X$, $d(\bar a,\bar b)= \lim_\cU d(a_i,b_i)$.  It is easy to see that $d$ is a pseudometric on $X$.  We define the \emph{ultraproduct of the family $(X_i,d_i)$ with respect to $\cU$}\index{ultraproduct}, denoted $\prod_\cU X_i$, to be the metric space obtained from $X$ by quotienting by the pseudometric $d$.  For $\bar x \in X$, we let $[\bar x]_\cU$ denote the equivalence class of $\bar x$ in the quotient. If $X_i$ is some fixed $X$ for all $i\in I$, we write $X^\cU$ for the ultraproduct and  call it the \emph{ultrapower}\index{ultrapower} of $X$ with respect to $\cU$.  There is a natural isometric embedding of $X$ into $X^\cU$, called the \emph{diagonal embedding}\index{diagonal embedding}, which sends $x \in X$ to the equivalence class of the constant sequence $(x : i \in I)$.

\begin{exercise}

\

\begin{enumerate}
    \item  Suppose that $(X_i:i\in I)$ is a collection of \emph{complete}, uniformly bounded metric spaces and $\cU$ is an ultrafilter on $I$.  Prove that $\prod_\cU X_i$ is complete. \item Suppose that $(X_n:n \in \bbN)$ is a collection of uniformly bounded metric spaces (which are not necessarily assumed to be complete) and $\cU$ is a non-principal ultrafilter on $\bbN$. Prove that $\prod_\cU X_n$ is complete.
    \end{enumerate}
\end{exercise}

\begin{exercise}
Suppose that $X$ is a compact bounded metric space.  Show that for any index set $I$ and ultrafilter $\cU$ on $I$, the diagonal embedding of $X$ into $X^\cU$ is surjective (and thus an isometry of metric spaces).
\end{exercise}

Throughout this article, we will be dealing with various metric spaces equipped with additional structure.  We highlight here a key point that will come up when we have to deal with such additional structure.  To keep things simple, suppose we have a complete, bounded metric space $(X,d)$  and, in addition, we have a unary function $f$ on $X$.  We can naturally define the product $X^\bbN$ with the $\sup$-metric and view $f$ as being defined co-ordinatewise on $X^\bbN$.  We ask:  what constraints must there be on $f$ if we wish to form an ultrapower of $(X,d,f)$ as a quotient of the product structure?  Notice in the classical case, there are no restrictions (one can always form the discrete ultrapower).  We claim that it is necessary and sufficient that $f$ be uniformly continuous.  Sufficiency is straightforward and is left as an exercise to the reader.  To verify necessity, suppose that there is $\e > 0$ such that, for every $n \in \bbN$ (see the footnote\footnote{We will adopt the convention that when $n \in \bbN$ is being treated as a number and not an index then $n \neq 0$.  This will avoid issues with division by 0.} on the use of $\bbN$ in this article), one can find $a_n,b_n \in X$ such that $d(a_n,b_n) < 1/n$ and yet $d(f(a_n),f(b_n)) \geq \e$.  Set $\bar a = (a_n : n \in \bbN)$ and $\bar b = (b_n : n \in \bbN)$.  If $\cU$ is a non-principal ultrafilter on $\bbN$, then $\lim_\cU d(a_n,b_n)=0$ and so $[\bar a]_\cU=[\bar b]_\cU$ as elements of the ultrapower $X^\cU$.  However, $\lim_\cU d(f(a_n),f(b_n))\geq \epsilon$ and so $[f(\bar a)]_\cU\not=[f(\bar b)]_\cU$, where $f(\bar a)=(f(a_n) \ : \ n\in \bbN)$ and likewise for $f(\bar b)$.  It follows that $f$ cannot be defined on the ultrapower as a quotient of the product structure.  

The conclusion we draw from the previous paragraph is that if we want to extend an ultraproduct construction to additional structure on a metric space, we are going to have to assume some form of uniform continuity.  In fact, a similar argument to the one just given shows that if one wants to construct an ultraproduct of uniformly bounded metric spaces $(X_n,d_n)$ for $n \in \bbN$ each of which is equipped with a uniformly continuous unary function $f_n$, the uniform continuity must be uniform along the sequence of structures, that is, for every $\e > 0$, there is $\delta >0$ so that for any $n\in \bbN$ and any $a,b\in X_n$, if $d_n(a,b) < \delta$, then $d_n(f_n(a),f_n(b)) \leq \e$ (see the discussion at the beginning of Section \ref{met-str} where we address the asymmetry in the above inequalities).

\section{Metric structures and languages}\label{met-str}

\subsection{Continuous languages and their interpretation}
We now introduce the structures we are interested in studying and the languages of the intended logic at the same time. A \emph{metric structure}\index{metric structure} $\M$, the continuous analog of a classical structure, will be multi-sorted.  The language will have a collection of sorts $\mathfrak{S}$.  For each sort $S \in \mathfrak{S}$, the language specifies a positive real number $B_S$.  The intended interpretation of $S$ in a metric structure $\M$ is a bounded, complete, non-empty metric space $S^\M$ whose diameter is at most $B_S$.  There is also a relation symbol $d_S$ in the language which is to be interpreted as the metric on the space $S^\M$.

A language also contains a set of function symbols $\mathfrak{F}$ and relation symbols $\mathfrak{R}$ (these are sometimes called predicate symbols in other presentations).  Each of these symbols comes with some additional information; we first treat function symbols.  If $f \in \mathfrak{F}$, then $f$ has a domain and co-domain as specified by the language.  The domain of $f$, denoted $dom(f)$, is a product of sorts $S_1 \times \cdots \times S_{n_f}$ for some $n_f \in \bbN$, and the co-domain of $f$, denoted $cod(f)$, is a sort $S_f$.  In addition, the language specifies an increasing, continuous function
 $\delta_f:[0,1] \rightarrow [0,1]$ satisfying $\lim_{\e \rightarrow 0^+} \delta_f(\e) = 0$; $\delta_f$ is  called the \emph{uniform continuity modulus}\index{uniform continuity modulus} of $f$.  The function $\delta_f$ is specified by the language (and motivated by the discussion in the previous section). 

\begin{comment}
    It may seem strong to insist that the uniform continuity modulus above is a continuous function.  As it turns out, it always can be arranged: 
\begin{prop}\cite[Proposition 2.10 and Remark 2.12]{BBHU}\label{2.10}
Let $F, G \colon X \rightarrow [0,1]$ be such
that
\[
(\forall \e > 0)( \exists \delta > 0)( \forall x \in X) (F(x) < \delta \text{ implies } G(x) \leq \e).
\]
Then there exists an increasing, continuous function $\alpha \colon [0, 1] \rightarrow [0, 1]$
such that $\alpha(0) = 0$ and such that,
for all $x \in X$, we have $G(x) \leq \alpha(F(x))$.  In fact, one may choose $\alpha$ so that it only depends on the map $\epsilon\mapsto \delta$ and not on the particular functions $F$ and $G$.
\end{prop}

Proposition \ref{2.10} is quite useful throughout continuous model theory.  For example, it is useful in expressing certain implications (which are otherwise difficult to express given that continuous logic is a ``positive'' logic); see \cite[Section 7]{BBHU}.  
We shall have occasion to use this proposition in Section \ref{sec-def} below.
\end{comment}
 
 In a metric structure $\M$ for this language, the function symbol $f$ is interpreted as a uniformly continuous function
\[
f^\M:S_1^\M \times \ldots \times S_{n_f}^\M \rightarrow S_f^\M
\] 
where the uniform continuity of $f^\M$ is witnessed by $\delta_f$, that is, for all $\e > 0$ and $a_i,b_i \in S_i^\M$  satisfying $d_{S_i}(a_i,b_i) < \delta_f(\e)$ for all $i=1,\ldots,n$, we have 
\[
d_{S_f}(f^\M(a_1,\ldots,a_{n_f}),f^\M(b_1,\ldots,b_{n_f})) \leq \e.
\]
The careful reader will notice the asymmetry between strict inequality and weak inequality in this definition of uniform continuity.  The short explanation for the asymmetry is that we want this definition to pass easily through the ultraproduct construction.  Nevertheless, the definition here is equivalent to one where we have strict inequalities on both sides.

Constant symbols will be thought of as function symbols with an empty sequence of sorts as domain and some sort as co-domain.  The interpretation in a metric structure will be a single element in the chosen co-domain.

For the relation symbols, if $R \in \mathfrak{R}$, the language specifies a domain, denoted $dom(R)$, which is a product of sorts $S_1\times\cdots \times S_{n_R}$ for some $n_R \in \bbN$.  The language also provides a bound $B_R > 0$ and a uniform continuity modulus $\delta_R$.  In a metric structure $\M$ for this language, $R$ is interpreted as a uniformly continuous function $R^\M:S_1^\M \times \ldots \times S_{n_R}^\M \to [-B_R,B_R]\subseteq \mathbb R$ with uniform continuity witnessed by $\delta_R$ as above.

To summarize, a metric language $\cL$ will contain:
\begin{itemize}
\item a set of sorts $\mathfrak{S}$ together with, for every $S \in \mathfrak{S}$, a relation symbol $d_S$ and number $B_S >0$,
\item a set of function symbols $\mathfrak{F}$ together with, for every $f \in \mathfrak{F}$, a product of sorts $dom(f)$, a sort $cod(f)$, and a uniform continuity modulus $\delta_f$, and
\item a set of relation symbols $\mathfrak{R}$ together with, for every $R \in \mathfrak{R}$, a product of sorts $dom(R)$, a number $B_R > 0$, and a uniform continuity modulus $\delta_R$.
\end{itemize}

A metric structure $\M$ for $\cL$ (also called an $\cL$-structure) specifies the following:
\begin{itemize}
\item for every sort $S \in \mathfrak{S}$, a complete non-empty metric space $S^\M$, with metric $d_S^\M$ of diameter bounded by $B_S$,
\item for every $f \in \mathfrak{F}$, a uniformly continuous function $f^\M$ whose uniform continuity is witnessed by $\delta_f$ and with domain $S_1^\M \times \cdots \times S_n^\M$ and co-domain $S^\M$, where $dom(f) = S_1 \times \ldots \times S_n$ and $cod(f)=S$, and
\item for every $R \in \mathfrak{R}$, a real-valued, bounded, uniformly continuous function $R^\M$ with uniform continuity modulus $\delta_R$, domain $S_1^\M \times \cdots \times S_n^\M$, where $dom(R) = S_1 \times \ldots \times S_n$, and co-domain $[-B_R,B_R]$.
\end{itemize}

For $\cL$-structures $\M$ and $\N$, we say that $\N$ is a \emph{substructure}\index{substructure} of $\M$ if
\begin{itemize}
\item for every $S \in \mathfrak{S}$, $S^\N$ is a closed subspace of $S^\M$,
\item for every $f \in \mathfrak{F}$, $f^\N$ is the restriction of $f^\M$ to $S_1^\N \times \ldots \times S_n^\N$, where $dom(f) = S_1 \times \ldots \times S_n$, and
\item for every $R \in \mathfrak{R}$, $R^\N$ is the restriction of $R^\M$ to $S_1^\N \times \ldots \times S_n^\N$, where $dom(R) = S_1 \times \ldots \times S_n$.
\end{itemize}

If $\cL' \subseteq \cL$ are two languages and $\M$ is an $\cL$-structure, then by $\M \restriction \cL'$ we mean the $\cL'$-structure obtained from $\M$ by only interpreting those sorts, relations and functions arising from $\cL'$.  We call $\M \restriction \cL'$ the \emph{reduct}\index{reduct} of $\M$ to $\cL'$ and we refer to $\M$ as an \emph{expansion}\index{expansion} of $\M \restriction \cL'$ to $\cL$.

We want to consider several categories of $\cL$-structures.  All of these categories will have $\cL$-structures themselves as objects; the difference between these categories will be the morphisms.  For now, we will consider two choices for maps between two $\cL$-structures, $\rho:\M \rightarrow \N$,  by which we mean a family of maps $(\rho_S : S \in \mathfrak{S})$ such that $\rho_S:S^\M \rightarrow S^\N$.  Below, we will write $\bar a \in \M$ to mean that $\bar a$ is a tuple from the union of the interpretations of the sorts of $\M$ and by $\rho(\bar a)$ we will mean $\rho_S^\M$ applied to the elements of the tuple $\bar a$ for the appropriate $S$'s.
\begin{itemize}
\item $\rho$ is a \emph{homomorphism}\index{homomorphism} if:
\begin{itemize}
    \item for all function symbols $f$ of $\cL$ and appropriately sorted $\bar a \in \M$, we have $\rho(f^\M(\bar a)) = f^\N(\rho(\bar  a))$;  
    \item for all relation symbols $R$ in $\cL$ and appropriately sorted $\bar a \in \M$, we have $R^\M(\bar a) \geq R^\N(\rho(\bar a))$.
    \end{itemize}
\item $\rho$ is an \emph{embedding}\index{embedding} if it is a homomorphism which furthermore satisfies, for all relation symbols $R$ in $\cL$ and appropriately sorted $\bar a \in \M$, $R^\M(\bar a) = R^\N(\rho(\bar a))$.  Notice that if $\rho$ is an embedding, then $\rho_S$ is an isometric embedding for all sorts $S$ of $\cL$.
\end{itemize}
In all categories we consider, $\rho$ is an \emph{isomorphism}\index{isomorphism} if it is an embedding and $\rho_S$ is surjective for all $S$.
If $\rho:\M \rightarrow \N$ is an embedding, then by $\rho(\M)$ we will mean the substructure of $\N$ induced on the image of the maps $\rho_S$ for all sorts $S$.

\subsection{Examples}

\begin{example}
A toy example of a metric structure is a bounded complete metric space.  The language has one sort and one binary relation symbol for the metric.  There are no uniform continuity requirements here.

A slightly more involved situation arises if one wants to consider an \emph{unbounded} complete metric space $(X,d)$.  While there is no canonical way to deal with this situation, the following approach is helpful, particularly when there is more structure as we will see later.  Fix a basepoint $\star \in X$.  Introduce in the language sorts $S_n$ for $n \in \bbN$.  $\{ a \in X : d(a,\star) \leq n \}$ is the intended interpretation of the sort $S_n$.  The metric on $S_n$ will just be the restriction of $d$ to this sort.  For $m<n$, we also introduce function symbols $i_{m,n}$ with domain $S_m$, co-domain $S_n$, and the identity function for modulus of uniform continuity.  The intended interpretation of $i_{m,n}$ is the natural inclusion of $S_m$ into $S_n$.
\end{example}

\begin{example}\label{ury-introduced}
A particularly interesting bounded metric space from the model-theoretic perspective is the Urysohn sphere $U$.  In \cite{U}, Urysohn constructed a separable complete $1$-bounded metric space which is both \emph{universal}, that is, every other separable 1-bounded metric space isometrically embeds into $U$, and \emph{ultrahomogeneous}, that is, if $A,B \subset U$ are finite and $f:A \rightarrow B$ is an isometric bijection, then $f$ can be extended to an automorphism of $U$.  $U$ is the unique  separable complete $1$-bounded metric space which is both universal and ultrahomogeneous.  Since $U$ is bounded, it is automatically an example of a metric structure.  

One way to build $U$ is via a Fra\"iss\'e-like construction using the class $\cal C$ of finite metric spaces with rational distances at most 1. (To see how Fra\"iss\'e constructions are done in the continuous setting, see Vignati's article in this volume.) There are only countably many objects in $\cal C$ and so, besides bookkeeping, one is led to the following problem:  given $A,B,C \in \cal C$ with $A\subset B$ and $A\subset C$, how can one amalgamate $B$ and $C$ over $A$?  To do this, one needs to decide the distance $d(b,c)$ for $b \in B\setminus A$ and $c \in C\setminus A$. We record the following lemma, whose proof is left as an exercise. 
\begin{lem} \label{met-amal}
Suppose that $A,B,C$ are all finite metric spaces with $A\subset B$ and $A\subset C$ and $B$ and $C$ are extensions of $A$.  Then $B$ and $C$ can be amalgamated over $A$ to form the metric space $B \cup_A C$ with the following properties:
\begin{itemize}
    \item the underlying set is $B\cup C$,
    \item the metric extends the metric on both $B$ and $C$, and
    \item the metric further satisfies 
\[
d(b,c) := \min_{a\in A}(d(b,a) + d(a,c))
\]
for $b \in B\setminus A$ and $c \in C\setminus A$.
\end{itemize}

\end{lem}

\begin{exercise}
Prove Lemma \ref{met-amal}.
\end{exercise}

Using Lemma \ref{met-amal}, one can form a chain of finite metric spaces from $\cal C$ in the style of a Fra\"iss\'e construction which yields a countable metric space $U_0$ with the property that whenever $A,B\in \cal C$ are such that $A \subset B$ and $A \subset U_0$, then there is a $B_0 \subset U_0$ with $A \subset B_0$ and an isometry $i:B \rightarrow B_0$ which is the identity on $A$. A mildly tedious exercise shows that the completion $U$ of $U_0$ satisfies the defining properties of the Urysohn sphere.
\end{example}

\begin{exercise}
Verify the claim made at the end of the previous example.
\end{exercise}

\begin{example}\label{CFOL}
A classical (discrete) structure can be thought of as a metric structure as follows.  Suppose that $\M$ is a classical structure in a classical language $\cL$. Equip each sort of $\M$ with the discrete metric (that is, $d(a,b) = 1$ if $a\neq  b$).  Additionally, for each relation symbol $R \in \cL$, reinterpret $R$ by the function $g_R = 1 - char(R^M)$, where $char(R^M)$ is the characteristic function of $R^M$. Note that $\bar a \in R^M$ iff $g_R(\bar a) = 0$.\footnote{Although arbitrary, the convention in continuous logic is that 0 represents true; see section 4 for further discussion.}  Moreover, $g_R$ is uniformly continuous as is the interpretation of every function symbol. For completeness, we should specify that the identity function is the modulus of continuity for all function and relation symbols. In this way, one can think of continuous logic as an extension of classical logic.
\end{example}

\begin{example}
Recall that a Hilbert space is a complex inner product space which is complete with respect to the induced metric.  We wish to capture the class of Hilbert spaces as an example of a class of metric structures. Here is the language $\cL_{HS}$ and standard interpretation we will work with:
\begin{enumerate}
\item For each $n\in \bbN$, we have a sort $S_n$.  The intended interpretation of $S_n$ in a Hilbert space is the ball centered at $0$ of radius $n$.  On $S_n$, there is a metric symbol with bound $2n$ whose intended interpretation is the metric induced by the inner product.
\item For each $n\in \bbN$, we have a binary function symbol $+_n$ with domain $S_n \times S_n$ and co-domain $S_{2n}$.  The intended interpretation for a Hilbert space is addition on the ball of radius $n$.  The uniform continuity modulus is $\e \mapsto \e/2$.
\item For each $n\in \bbN$, we have a constant symbol $0_n$ in sort $S_n$ which will be interpreted as 0.  
\item For each $n \in \bbN$ and $\lambda \in \bbC$, there is a unary function symbol $\lambda_n$ with domain $S_n$ and co-domain $S_m$, where $m = \lceil |\lambda| \rceil$. The uniform continuity modulus can be taken to be $\e \mapsto \e/m$  The intended interpretation is scalar multiplication by $\lambda$.
\item For each $m,n\in \bbN$ with $m < n$, we have unary function symbols $i_{m,n}$ with domain $S_m$ and co-domain $S_n$.  The intended interpretation is the inclusion map from the ball of radius $m$ into the ball of radius $n$ and the uniform continuity modulus is the identity function.
\item For each $n \in \bbN$, we have two binary relation symbols $Re(\langle \cdot,\cdot \rangle)_n$ and $Im(\langle \cdot,\cdot \rangle)_n$ on $S_n$.  The intended interpretation of these symbols is the restriction to $S_n$ of the real and imaginary parts of the inner product.  The bounds for each of these symbols are $n^2$ and the moduli of uniform continuity can be take to be $\e \mapsto \e/2n$.
\end{enumerate}
If $H$ is a Hilbert space, then the $\cL_{HS}$-structure just described will be denoted by $D(H)$ and referred to as the \emph{dissection of $H$}\index{dissection}.
In practice, we will write formulas in this language which will omit the subscripts when the intention is clear.  We may also write the inner product symbols without reference to the real and imaginary part if it is clear that an equivalent formula can be written adhering to the above formalism.  

\end{example}

One of the most important examples for the current volume is the class of \cstar-algebras.  Let us look at this class in detail.
\begin{example}\label{C*ex}
Given a \cstar-algebra $A$, the natural metric on $A$ is that induced by the operator norm.  As in the case of Hilbert spaces treated in the previous example, this metric is unbounded.  Once again, we get around this issue by introducing a series of sorts to represent operator norm balls of various radii.  As before, we fix sorts $S_n$ for each $n \in \bbN$ with the intended interpretation being the operator norm ball of radius $n$ in $A$.  The intended metric on each sort is just the restriction of the metric induced by the operator norm to the given sort.  The language will also have a number of function symbols.  Since we will be dissecting  \cstar-algebras, as in the case of Hilbert spaces, we will need to consider the restriction of functions like addition and multiplication to each sort and we need to connect the sorts via embedding maps.  With that in mind, we introduce the language $\cL_{C^*}$ with its intended interpretation:
\begin{enumerate}
\item For each $n \in \bbN$, there are binary function symbols $+_n$ and $\cdot_n$ with domain $S_n \times S_n$ and co-domains $S_{2n}$ and $S_{n^2}$ respectively.  The intended interpretation is the restriction of addition and multiplication on $A$ to the operator norm ball of radius $n$. The uniform continuity modulus for $+_n$ is $\epsilon \mapsto \epsilon/2$ and for $\cdot_n$ is $\epsilon \mapsto \epsilon/2n$.
\item For each $n \in \bbN$, there is a unary function symbol $\cdot^{\ast_n}$ with domain and co-domain $S_n$.  The intended interpretation is the restriction of the adjoint to the operator norm ball of radius $n$.  The uniform continuity modulus is the identity function.
\item For each $n \in \bbN$, there is a constant symbol $0_n$ in $S_n$.   The intended interpretation is 0.
\item For each $\lambda \in \bbC$ and $n \in \bbN$, we have a unary function symbol $\lambda_n$ with domain $S_n$ and co-domain $S_m$, where $m = \lceil |\lambda | \rceil$.  The uniform continuity modulus of this symbol is $\epsilon \mapsto \epsilon/m$.
\item For each $n \in \bbN$, the metric symbol on $S_n$, denoted $d_n$, will be interpreted as $\|x - y\|$ for $x,y$ in the operator norm ball of $A$ of radius $n$.
\item For each $m,n\in \bbN$ with $m < n$, there is a unary function symbol $i_{m,n}$ with domain $S_m$ and co-domain $S_n$.  The intended interpretation is the inclusion map between $S_m^A$ and $S_n^A$.  The uniform continuity modulus is the identity map.

\end{enumerate}

Given a \cstar-algebra $A$, the $\cL_{C^*}$-structure considered above is denoted $D(A)$ and is called the \emph{dissection of $A$}.  It is clear that if we start with $A$ and consider $D(A)$, we can then recover $A$ from $D(A)$ (simply take the union of the sorts and piece together the operations as given).  Moreover, $D(A) \cong D(B)$ (as $\cL_{C^*}$-structures) if and only if $A \cong B$ (as \cstar-algebras).

We now address a subtle point.  Suppose that we have a \textbf{$\cL_{C^*}$}-substructure $\M \subseteq D(A)$ for some \cstar-algebra $A$.  Does $\M$ arise as $D(B)$ for some subalgebra $B \subseteq A$?  To be clear, it is easy to see that if one considers the direct limit of the sorts $S_n^\M$ along the embeddings arising from the maps $i_{m,n}$, we do indeed obtain a \cstar-subalgebra $B$ of $A$ (the completeness of each individual sort guarantees that the union is complete).  What is in question is whether each individual sort $S_n$ captures exactly the elements of $B$ which are of operator norm at most $n$.  In order to see that we can arrange for this to happen, we need the following fact. The functional analytic argument that underlies this fact can be found in \cite[section 3.1 at the top of page 486]{mtoa2}.

 \begin{fact}\label{func-analytic}
 For any $n\in \bbN$, there are one-variable polynomials $q^n_k$ ($k\in \bbN$) such that if $A$ is a \cstar-algebra and $a\in A$ satisfies $\|a\| \leq n$, then $\|q_k^n(a)\| \leq 1$ for all $k \in \bbN$.  If moreover $\|a\| \leq 1$, then $q_k^n(a)\to a$ as $k\to \infty$.
 \end{fact}  
 We henceforth assume that, for each $k,n\in \bb N$ we have added function symbols $\mathbf{q^n_k}$  to the language $\cL_{C^*}$ with domain $S_n$ and co-domain $S_1$. The modulus of continuity can be taken to be the modulus of continuity of the polynomial $q^k_n$ when restricted to the operator norm ball of radius $n$. In the dissection of a \cstar-algebra, we will interpret $\mathbf{q^n_k}$ by the value of the corresponding polynomial (taking values in $S_1$). Since each polynomial $q^n_k$ is naturally a composition of the functions in the language of \cstar-algebras, there is some $m \geq n$ such that $S_m$ will be the (na\"ive) co-domain of $q^n_k$ when restricted to $S_n$. When $\mathbf{q^n_k}$ is interpreted in the dissection of a \cstar-algebra, we have
 \[
 i_{1,m}(\mathbf{q^n_k}(x)) = q^n_k(x)
 \]
 for all $x$ in sort $S_n$.
 The justification for adding these function symbols is as follows.
 
 Suppose that $\M$ and $B$ are as above and $\| a \| \leq 1$ for some $a \in B$.  Take $n \in \bbN$ for which $a \in S_n^\M$.  Consider $q^n_k$ for $k \in \bbN$ as in Fact \ref{func-analytic}.  Then in $B$, $q^n_k(a)$ tends to $a$ in operator norm.  On the other hand, since the co-domain of $\mathbf{q^n_k}$ is $S_1$, $(\mathbf{q^n_k})^\M(a)$ converges to something in $S_1^\M$, call it $b$.  Since $i_{1,n}$ is an isometric embedding, $i_{1,n}^\M(\mathbf{q^n_k}(a))$ converges to $a$ in $S_n^\M$ and hence $i_{1,n}^M(b) =a$, which means that $a$ and $b$ are identified in $B$.  We conclude that $D(B) \cong M$.

 The dissection functor $D$ is in fact an equivalence of categories from the category of \cstar-algebras with homomorphisms as morphisms and $\cL_{C^*}$-structures of the form $D(A)$ for \cstar-algebras $A$, again with homomorphisms as morphisms.  The equivalence between these categories is actually stronger than this, but this stronger statement will have to wait until we define ultraproducts of structures.

\end{example}

\section{Formulas}\label{sec-form}

\subsection{Basic Formulas}
We now introduce the syntax of formulas for continuous logic.  Fix a continuous language $\cL$ with sorts $\mathfrak{S}$, function symbols $\mathfrak{F}$, and relation symbols $\mathfrak{R}$.

For each sort $S \in \mathfrak{S}$, we introduce variables $x_n^S$ for $n \in \bbN$.  As in classical logic, we need to define terms before formulas.
{\em Terms} in $\cL$, along with their domain, co-domain, and associated uniform continuity modulus, are defined inductively as follows:
\begin{enumerate}
\item If $x$ is a variable of sort $S$, then $x$ is a term.  Its domain and co-domain are $S$ and its uniform continuity modulus is the identity function.
\item If $f \in \mathfrak{F}$ is a function symbol with domain $S_1 \times \ldots \times S_n$ and co-domain $S$ and $\tau_1,\ldots,\tau_n$ are terms with co-domains $S_1,\ldots, S_n$ respectively, then $f(\tau_1,\ldots,\tau_n)$ is a term with co-domain $S$ and domain determined by composition.  The uniform continuity modulus of this term is also determined by composition; for instance, one can take
\[
\epsilon \mapsto \min\{ \delta_1(\delta_f(\epsilon)/2),\ldots, \delta_n((\delta_f(\epsilon)/2) \},
\]
where $\delta_f$ is the uniform continuity modulus of $f$ and $\delta_k$ is the uniform continuity modulus of $\tau_k$ for every $k = 1, \ldots, n$.
\end{enumerate}

\begin{exercise}
Verify that the function displayed at the end of item (2) above is indeed a uniform continuity modulus for $f(\tau_1,\ldots,\tau_n)$.
\end{exercise}

The \emph{basic formulas}\index{basic formula} of $\cL$ are now defined by induction. We use the term basic formula because, as we will see later in this article, we will want to introduce a more general notion of formula to accommodate the approximate nature of continuous logic.\footnote{We hope there is no confusion with this use of the term basic formula and its various uses in classical logic.}  Each basic formula will have a domain, bound, and uniform continuity modulus as is the case for relation symbols.  
\begin{enumerate}
\item Suppose that $R \in \mathfrak{R}$ has domain $S_1\times\ldots\times S_n$ and bound $B_R$, and $\tau_1,\ldots,\tau_n$ are terms with co-domains $S_1,\ldots, S_n$ respectively.  Then $R(\tau_1,\ldots,\tau_n)$ is a basic formula with bound $B_R$, and domain and uniform continuity modulus determined by composition.  In particular, since all the symbols for metrics on the sorts are included in $\mathfrak{R}$, if $\tau_1$ and $\tau_2$ are sorts with the same co-domain $S$, then $d_S(\tau_1,\tau_2)$ is a basic formula.
\item If $f:\bbR^n \rightarrow \bbR$ is a continuous function and $\varphi_1, \ldots,\varphi_n$ are basic formulas, then $\psi := f(\varphi_1,\ldots,\varphi_n)$ is also a basic formula.  If $B_1,\ldots,B_n$ are the bounds of $\varphi_1,\ldots,\varphi_n$ respectively, then the bound of $\psi$ is the minimum $B$ such that 
\[
f(\prod_{i = 1}^n [-B_i,B_i]) \subseteq [-B,B].
\]
The domain of $\psi$ is determined by composition. The uniform continuity modulus of $\psi$ is also determined by composition from the uniform continuity moduli of the $\varphi_i$'s and a fixed uniform continuity modulus of $f$ restricted to $\prod_{i = 1}^n [-B_i,B_i]$.
\item Suppose that $x$ is a variable and $\varphi$ is a formula.  We can then form two {\em quantified} basic formulas $\sup_x \varphi$ and $\inf_x \varphi$.  The bound and uniform continuity modulus of these formulas is the same as $\varphi$.  The domain is only effected by the possible removal of the sort associated with $x$.
\end{enumerate}
The basic formulas arising from the first clause are called \emph{atomic formulas}\index{atomic formula}.  Basic formulas arising from the first two clauses are called \emph{quantifier-free}\index{quantifier-free formula} basic formulas. 

\begin{comment}
In the second clause above, we have chosen to be very generous about the functions - the connectives - we allow.  It is convenient theoretically to allow all continuous functions (which are actually uniformly continuous since we are only really interested in their values on some hypercube in $\bbR^n$).  It is sometimes worthwhile to work with a smaller collection of connectives.  The expressive power of continuous logic would not change if we used a collection of connectives which were dense in the set of all continuous functions.  For instance, we could restrict ourselves to using only polynomials with rational coefficients as connectives.  Other restricted sets of connectives have been suggested (see \cite[Section 6]{BBHU} for more details).  As we will see, it is the density character of the set of formulas which matters the most in continuous logic and not the cardinality of the set of formulas.

Another set of basic formulas to highlight are the \emph{prenex}\index{prenex formula} formulas, which are those basic formulas beginning with a block of quantifiers followed by a quan-tifier-free basic formula. In a natural (pseudo)metric defined on the space of all basic formulas (see Subsection \ref{sec-mt1}), the prenex formulas are dense.  We will come back to prenex formulas in Subsection \ref{SW}.
\end{comment}

We now consider the interpretation of basic formulas in a metric structure.  Fix an $\cL$-structure $\M$.

A term $\tau$ is interpreted in $\M$ by $\tau^\M$ exactly as it would be in classical logic by following the inductive definition and composing.

\begin{exercise}
    Verify that the domain, co-domain and uniform continuity modulus of $\tau$ as known to the language are in fact the domain, co-domain and uniform continuity modulus of $\tau^\M$.
\end{exercise}

The interpretation of basic formulas is not much more complicated.  Suppose that $\bar a = (a_1,\ldots,a_n)$ is a sequence of sorted elements of $\M$ and the free variables $\bar x = (x_1,\ldots, x_n)$ of a basic formula $\varphi$ are similarly sorted.  We wish to define the interpretation $\varphi^\M(\bar a)$ inductively:

\begin{enumerate}
\item If $\varphi := R(\tau_1,\ldots,\tau_k)$ for terms $\tau_1,\ldots,\tau_k$, then
\[ 
\varphi^\M(\bar a) := R^\M(\tau_1^\M(\bar a),\ldots,\tau_k^\M(\bar a)).
\]

\item If $\varphi := f(\psi_1,\ldots,\psi_k)$ for basic formulas $\psi_1,\ldots,\psi_k$ and a continuous real-valued function $f$, then 
\[
\varphi^\M(\bar a) := f(\psi^\M_1(\bar a),\ldots,\psi^\M_k(\bar a)).
\]
\item If $\varphi := \sup_x \psi$, then 
\[
\varphi^\M(\bar a) := \sup\{ \psi^\M(b,\bar a) \colon b \in S^\M \}
\]
where $S$ is the sort of $x$.  Notice that this is well-defined since the co-domain of $\psi^\M$ is bounded.

\noindent Similarly, if $\varphi := \inf_x \psi$, then
\[
\varphi^\M(\bar a) := \inf\{ \psi^\M(b,\bar a) \colon b \in S^\M \}
\]
where again $S$ is the sort of $x$.  

\
\end{enumerate}

\begin{exercise}
    Verify that in each case of the above definition, these formulas have bound, domain and uniform continuity modulus known to the language.
\end{exercise}

\subsection{Some examples of formulas}
Since it takes some getting used to writing expressions in continuous logic, let us look at some standard ways of expressing certain statements.
\begin{example}
Suppose we have two terms $\tau(\bar x)$ and $\sigma(\bar x)$ with the same co-domain.  How do we assert that they are equal when they are interpreted?  If $d$ is the metric symbol on the sort of the common co-domain of the terms, then the sentence $\sup_{\bar x} d(\tau(\bar x),\sigma(\bar x))$ evaluates to $0$ in some metric structure $\M$ if and only if the functions $\tau^\M$ and $\sigma^\M$ are equal.
\end{example}

\begin{example}
We sometimes wish to assert that one formula is less than or equal to another in value.  We do this with the help of the truncated substraction function $x \dotminus y := \max\{0,x-y\}$.  Note that $\dotminus$ is a continuous function from $\bb R^2$ to $\bb R$, hence may be used as a connective.  Suppose that $\varphi(\bar x)$ and $\psi(\bar x)$ are two formulas.  Then the sentence
\[
\sup_{\bar x} (\varphi(\bar x) \dotminus \psi(\bar x))
\] evaluates to $0$ in a structure $\M$ if and only if $\varphi^\M\leq \psi^\M$ as functions.
\end{example}

\begin{example}\label{ex4.3}
Another useful thing to be able to write out is the isomorphism type of some finite metric space.  Suppose that $A = \{a_1,\ldots,a_n\}$ is a finite metric space with $n$ distinct elements and metric $\delta$.  Set $r_{ij} = \delta(a_i,a_j)$ for $1 \leq i < j \leq n$.  Consider the formula 
\[
D_A(x_1,\ldots,x_n) = \max_{1 \leq i < j \leq n} |d(x_i,x_j) - r_{ij}|,
\]
where $x_1,\ldots,x_n$ are variables of sort $S$ and $d$ is the metric symbol for that sort.  Note that $\max$ is a continuous function (here of $n(n-1)$ inputs) as is $a \mapsto |a - r|$, hence they may be used as connectives in the formation of the formula $D_A$ above. If we have a structure $\M$ with sort $S$ and $b_1,\ldots,b_n \in S^\M$, then $\varphi^\M(b_1,\ldots,b_n) =0$ if and only if the map $a_i \mapsto b_i$ for $i =1,\ldots,n$ is an isometric embedding of $A$ into $S^\M$.
\end{example}

\section{Basic Continuous Model Theory}\label{sec-mt}

\subsection{Some terminology}\label{sec-mt1}

Before defining a more general notion of formula, we need to introduce some logical terminology which roughly parallels the classical case.

Fix a metric language $\cL$. A {\em sentence} in $\cL$ is an $\cL$-formula with no free variables.  The set of $\cL$-sentences will be denoted $Sent(\cL)$. A set of sentences $\Gamma$ is {\em satisfiable} if there is an $\cL$-structure $\M$ so that $\varphi^\M = 0$  \footnote{The choice of 0 here is somewhat arbitrary.  If we were interested in whether a certain sentence $\varphi$ evaluated to some real number $r$ then, we could equivalently ask if the sentence $\varphi - r$ was satisfied. }  for all $\varphi \in \Gamma$, in which case we say that $\M$ \emph{satisfies $\Gamma$}\index{satisfaction} and we write $\M \models \Gamma$.  When $\Gamma=\{\varphi\}$, we write $\M\models \varphi$ instead of $\M\models \{\varphi\}$ and say that $\M$ satisfies $\varphi$. We extend this notation by writing, for $\Gamma$ a set of sentences and $\varphi$ a sentence, $\Gamma \models \varphi$ to mean that whenever $\M\models \Gamma$, then $\M\models \varphi$.

A set of satisfiable ($\cL$-)sentences $T$ is is called an \emph{($\cL$)-theory}\index{theory}. We say that a class of $\cL$-structures $\cal C$ is \emph{elementary}\index{elementary} (or \emph{axiomatizable}\index{axiomatizable}) if there is an $\cL$-theory $T$ such that $\cal C$ is the collection of models of $T$. 

For an $\cL$-structure $\M$, 
%the \emph{theory of $\M$}, denoted $Th(\M)$, is the set of all $\cL$-sentences $\varphi$ such that $\varphi^\M = 0$.  
consider the linear \emph{evaluation map}
\[
ev_\M \colon Sent(\cL) \rightarrow \bbR \text{ defined by } ev_\M(\varphi) = \varphi^\M.
\]
Notice that if $ev_M(\varphi) = r$, then $ev_M(\varphi - r) = 0$, hence $ev_\M$ is completely determined by its kernel.  With this in mind, it is convenient to define \emph{the theory of $\M$}, $Th(\M)$, as the set of sentences $\varphi$ such that $\varphi^\M = 0$.   Notice also that $ev_\M$ is more than linear: if $\varphi_1,\ldots,\varphi_n$ are $\cL$-sentences and $f$ is a continuous real-valued function on $\bbR^n$, then
\[
ev_\M(f(\varphi_1,\ldots,\varphi_n)) = f(ev_\M(\varphi_1),\ldots,ev_\M(\varphi_n)).
\]

Two $\cL$-structures $\M$ and $\N$ are \emph{elementarily equivalent}\index{elementary equivalent}, written $\M \equiv \N$, if $Th(\M) = Th(\N)$.  A set of sentences $\Gamma$ is \emph{complete}\index{complete theory} if it is satisfiable, and whenever $\M$ and $\N$ satisfy $\Gamma$, then $\M \equiv \N$. Notice that for an $\cL$-structure $\M$, $
Th(\M)$ is a complete ($\cL$-)theory.  If $\N$ is a substructure of $\M$, then we say that $\N$ is an \emph{elementary substructure}\index{elementary substructure} of $\M$, written $\N \preceq \M$, if whenever $\varphi(\bar x)$ is an $\cL$-formula and $\bar a \in \N$ is a sequence sorted in the same manner as the variables $\bar x$, then $\varphi^\N(\bar a) = \varphi^\M(\bar a)$.  We write $\N \prec \M$ if $\N\preceq \M$ and $\N \neq \M$. If $\rho:\M \rightarrow \N$ is an embedding, then it is an \emph{elementary embedding}\index{elementary embedding} if $\rho(\M) \preceq \N$.

We can now define another category of $\cL$-structures.  Given an $\cL$-theory $T$, $Mod(T)$ is the category with models of $T$ as objects and elementary embeddings as morphisms.  If the theory $T$ is empty, that is, the collection of objects is just the collection of all $\cL$-structures, then we refer to this category as $Str(\cL)$.

\subsection{General Notion of Formula}

We now generalize the notion of basic formula.  It is important to note that this general notion of formula is relative to a background theory.  The general notion of formula follows from introducing the following semi-norm on the set of basic formulas.  Fix an $\cL$-theory $T$ and a fixed set of variables $\bar x$.  Given a basic formula $\varphi(\bar x)$, define $\| \varphi(\bar x) \|_T = \sup \{ |\varphi^\M(\bar a)| \colon \bar a \in \M, \M \models T \}$.  Let $Form_{\bar x}(T,\cL)$ be the completion of the set of basic $\cL$-formulas in free variables $\bar x$ with respect to the pseudometric induced by the semi-norm $\| \cdot \|_T$ (that is, the distance between $\varphi$ and $\psi$ is $\|\varphi-\psi\|_T$).  We say $\varphi$ and $\psi$ are \emph{$T$-equivalent}\index{$T$-equivalent} if $\|\varphi-\psi\|_T = 0$.  Thus, $Form_{\bar x }(T)$ is obtained by first identifying $T$-equivalent formulae and then completing the resulting metric space.  Elements of $Form_{\bar x}(T,\cL)$ will be called \emph{$T$-formulas}\index{$T$-formula} and are limits of basic formulas in the sense of the semi-norm $\| \cdot \|_T$.  If $T = \emptyset$, we will refer to elements of $Form(\emptyset,\cL)$ simply as \emph{$\cL$-formulas}\index{$\cL$-formula}. If $\varphi$ is a $T$-formula which is the limit of basic formulas $\varphi_n$, then we can interpret $\varphi$ in a model $\M$ of $T$ by setting $\varphi^\M = \lim_n \varphi_n^\M$.

\begin{exercise}
Verify that the definition of $\varphi^\M$ for $\varphi$ a $T$-formula given above is independent of the choice of representative sequence.  Moreover, prove that for each $T$-formula $\varphi$, there is $B>0$ such that $\varphi^\M$ takes values in $[-B,B]$ for all models $\M$ of $T$.
\end{exercise}

Notice that for a theory $T$, the collection of $T$-formulas is closed under composition with continuous functions i.e. if $\varphi_1,\ldots,\varphi_n$ are $T$-formulas and $f:\bbR^n \rightarrow \bbR$ is a continuous function then $f(\varphi_1,\ldots,\varphi_n)$ is a $T$-formula. Moreover, the set of $T$-formulas is closed under $\sup$ and $\inf$.  That is, if $\varphi(x,\bar y)$ is a $T$-formula then so is $\sup_x \varphi(x,\bar y)$ and $\inf_x \varphi(x,\bar y)$.

\begin{example}\label{weighted}
    Suppose that we have a theory $T$ and $T$-formulas $\varphi_n$ in variables $\bar x$ with bounds $B_n$ for $n \in \bbN$.  Then the following is a $T$-formula in the sense just defined:
    \[
    \sum_{n \in \bbN} \frac{\varphi(\bar x)}{B_n2^n}.
    \]
    It is not hard to see that all $T$-formulas can be realized as weighted sums of basic formulas.
\end{example}

$Form_{\bar x}(T,\cL)$ has the structure of a real normed algebra: it contains the constant functions and is closed under addition and multiplication.  
It will be convenient to generalize this notion as follows. 
\begin{defn}\label{real-alg}
A \emph{real algebra system $A$ of $T$-formulas}\index{real algebra system} consists of, for every sequence of variables $\bar x$, a closed subalgebra $A_{\bar x}$ of $Form_{\bar x}(T,\cL)$.
\end{defn}  

We have three particular cases of real algebra systems in mind.
\begin{example}

\

\begin{enumerate}
    \item Suppose that $\cL' \subseteq \cL$.  Given a sequence $\bar x$ of variables, define $A_{\bar x}$ to be the closure in $Form_{\bar x}(T,\cL)$ of the set of basic $\cL'$-formulas. $A$ is then a real algebra system of $T$-formulas.
    \item Given a sequence $\bar x$ of variables, let $A_{\bar x}$ be the closure in $Form_{\bar x}(\emptyset,\cL)$ of the set of quantifier-free basic $\cL$-formulas. Then $A$ is a real algebra system of $\cL$-formulas whose elements are called \emph{quantifier-free formulas}.
    \item Let $A_{\bar x}$ be the closure in $Form_{\bar x}(\emptyset,\cL)$ of the set of prenex formulas in $\cL$ in the free variables $\bar x$.  Then $A$ is a real algebra system.
    
\end{enumerate}
\end{example}

For a topological space $X$, we let $\chi(X)$ denote the \emph{density character of $X$}\index{density character}, namely the smallest cardinality of a dense subset of $X$.  We will refer to the \emph{density character of a language $\cL$ with respect to a theory $T$}, denoted by $\chi(T,\cL)$, as 
\[
\sum_{\bar x} \chi(Form_{\bar x}(T,\cL)),
\]
where the sup is over all tuples $\bar x$ of variables.
If $T$ is the empty theory, we will omit it from the notation and just write $\chi(\cL)$.  If $\chi(\cL) = \aleph_0$, we say that $\cL$ is \emph{separable}.  We also define the density character of an $\cL$-structure $\M$, denoted  $\chi(\M)$, to be
\[
\sup_{S \in \mathfrak{S}} \chi(S^\M).
\]

\subsection{First model-theoretic results}

\begin{prop}[Tarski-Vaught test]
Suppose that $\N \subseteq \M$.  Then the following are equivalent:
\begin{enumerate}
\item $\N \preceq \M$.
\item  For all basic $\cL$-formulas $\varphi(y,\bar x)$ and all $\bar a \in \N$,
\[
\inf\{ \varphi^\M(b,\bar a) \colon b \in \N \} = \inf\{ \varphi^\M(b,\bar a) \colon b \in \M \}.
\]
\end{enumerate}
\end{prop}

\proof (1) implies (2) by the definition of elementary substructure.  For (2) implies (1), we argue by induction on the formation of basic formulas.  The interesting case is that of an $\inf$ quantifier. Thus, suppose $\varphi(\bar x) = \inf_y \psi(y,\bar x)$ and $\bar a \in \N$.  By induction, we always have $\varphi^\N(\bar a) \geq \varphi^\M(\bar a)$.  To prove the reverse inequality, fix $\e > 0$ and take $b \in \M$ such that $\psi^\M(b,\bar a) < \varphi^\M(\bar a) + \e$.  By condition (2), we can find $b' \in \N$ such that $\psi^\M(b',\bar a) < \varphi^\M(\bar a) + \e$, and so by induction and letting $\e$ tend to 0, we have $\varphi^\N(\bar a) \leq \varphi^\M(\bar a)$. \qed

\begin{thm}[Downward L\" owenheim-Skolem]
Given $X \subseteq \M$, there is $\N \preceq \M$ such that $X \subseteq N$ and $\chi(\N)\leq \chi(X) + \chi(\cL)$.
\end{thm}

\proof We accomplish the creation of $\N$ by iteratively closing off in order to satisfy the Tarski-Vaught test; we only sketch of the details.  
Fix a set of formulas $\{ \varphi_i(y,\bar x_i) \colon i \in I\}$ which is dense in $Form(Th(\M),\cL):=\bigcup_{\bar x}Form_{\bar x}(Th(\M),\cL)$ with $|I| = \chi(\cL)$.  We create an increasing sequence $(X_n)$ of subsets of $\M$ such that:
\begin{enumerate}
\item $X_0 = X$,
\item for each $n \in \bbN$, $i \in I$, and $\bar a \in X_n$, there is a sequence $(b_n)$ contained in $X_{n+1}$ such that 
\[
\textstyle{(\inf_y \varphi(y,\bar a))^\M = \inf\{ \varphi^\M(b_n,\bar a) : n \in \bbN \}},
\]
and
\item $\chi(X_n) \leq \chi(X) + \chi(\cL)$.
\end{enumerate}
We let $N$ be the closure (in $M$) of the union of the $X_n$'s and leave it as an exercise to the reader to check that $N$ is the universe of a substructure $\N$ of $\M$ which is in fact an elementary substructure of $\M$. \qed

\begin{exercise}
Verify the claim made at the end of the previous proof.
\end{exercise}

\begin{defn}
Fix an $\cL$-structure $\M$ and a subset $A \subseteq \M$.  The language $\cL_A$ is $\cL$ together with a new constant symbol $c_a$ for every $a \in A$  (and where $c_a$ has the same sort as $a$).
$\M$ can be canonically expanded to an $\cL_A$-structure $\M_A$ by interpreting $c_a$ as $a$ for every $a \in A$.  

We define two types of $\cL_M$-theories (often referred to as diagrams):
\begin{enumerate}
\item The \emph{atomic diagram}\index{atomic diagram} of $\M$, denoted $Diag(\M)$,  is the set
\[
\{ \varphi \colon \varphi \text{ is an }  \text{ an atomic $\cL_\M$-sentence and } \M_\M \models \varphi \}.
\]
\item The \emph{elementary diagram}\index{elementary diagram} of $\M$, $Elem(\M)$, is the set 
\[
\{ \varphi \colon \varphi \text{ is an } \text{ an basic $\cL_\M$-sentence and } \M_\M \models \varphi \}.
\]
\end{enumerate}
\end{defn}

\begin{prop}
For $\cL$-structures $\M$ and $\N$ we have that:
\begin{enumerate}
\item $\M$ can be embedded into $\N$ if and only if $\N$ can be expanded to an $\cL_\M$-structure which satisfies the atomic diagram of $\M$.
\item $\M$ can be elementarily embedded into $\N$ if and only if $\N$ can be expanded to an $\cL_\M$-structure which satisfies the elementary diagram of $\M$.
\end{enumerate}
\end{prop}

\proof The proofs are similar so we will only prove (2). If $f\colon \M \rightarrow \N$ is an elementary embedding, then by interpreting $c_m$ in $\N$ by $f(m)$, it is straightforward to show this expansion of $\N$ satisfies the elementary diagram of $\M$.  In the other direction, if $\tilde{\N}$ is an $\cL_\M$ expansion of $\N$ which satisfies the elementary diagram of $\M$, then define $f:\M \rightarrow \N$ by $f(m) = c^\N_m$ for all $m \in \M$.  Again it is straightforward to show that this is an elementary map. \qed

The following result is analogous to the classical result regarding unions of chains.
\begin{prop}\label{direct-limit}
The category of $\cL$-structures with elementary embeddings as morphisms has directed colimits.  More precisely, given:
\begin{enumerate}
\item a directed partial order $(P,<)$ and a family of $\cL$-structures $(\M_i \colon i \in P)$, and
\item elementary maps $\rho_{ij} : \M_i \rightarrow \M_j$ for every $i < j$ satisfying 
$\rho_{jk}\circ \rho_{ij} = \rho_{ik}$ whenever $i < j < k$,
\end{enumerate}
 then there is an $\cL$-structure $\N$ and elementary maps $\nu_i : \M_i \rightarrow \N$ such that 
 \begin{enumerate}
 \item $\nu_i = \nu_j \circ \rho_{ij}$ for $i < j$, and 
 \item $\N$ is the closure of $\bigcup_{i\in P}\nu_i(\M_i)$.
 \end{enumerate}
\end{prop}

\begin{exercise}
    Prove the previous proposition.
\end{exercise}

\subsection{Ehrenfeucht-Fra\"iss\'e games}
We now introduce the continuous version of Ehrenfeucht-Fra\"iss\'e games which provide a way of recognizing elementary equivalence.  These games have been useful in the study of theories of various operator algebras (see for example \cite{GH-McDuff, GS}).

Fix a set of basic $\cL$-formulas $\Delta$ in the variables $\bar x = (x_i \colon 1 \leq i \leq n)$, as well as an $\e > 0$.  Suppose $\M$ and $\N$ are two $\cL$-structures.  We describe a game, the $EF(\M,\N,\Delta,\e)$-game,  between two players, Player I and Player II.  The game has $n$ rounds.  In the i$^{th}$ round, Player I chooses either $a_i \in \M$ or $b_i \in \N$ of whatever sort the variable $x_i$ belongs to and then Player II replies with $b_i \in \N$ or $a_i \in \M$ respectively.  After $n$ rounds, there are two sequences $\bar a = (a_1,\ldots,a_n)$ and $\bar b = (b_1,\ldots,b_n)$.  Player II wins the game if for every $\varphi(\bar x) \in \Delta$, we have
\[
|\varphi^\M(\bar a) - \varphi^\N(\bar b)| < \e.
\]

\begin{thm}\label{EF}
For two $\cL$-structures $\M$ and $\N$, the following are equivalent:
\begin{enumerate}
\item $\M \equiv \N$.
\item For every finite set $\Delta$ of atomic $\cL$-formulas and every $\e > 0$, Player II has a winning strategy for the $EF(\M,\N,\Delta,\e)$-game.
\item For every finite set $\Delta$ of basic $\cL$-formulas and every $\e > 0$, Player II has a winning strategy for the $EF(\M,\N,\Delta,\e)$-game.

\end{enumerate}
\end{thm}

\begin{exercise}
Prove Theorem \ref{EF}.  (Hint:  For the dirction (3) implies (1), proceed by induction on the natural notion of \emph{quantifier depth} of a basic $\cL$-formula.) 
\end{exercise}

A reader looking for a solution to the previous exercise can consult  \cite{H}.

\subsection{Examples of continuous theories}

\begin{example}[The Urysohn sphere]\label{ury-example}

%For a finite metric space $A = \{a_1,\ldots,a_n\}$ with metric $d$ let the formula $D_A(x_1,\ldots,x_n)$ be
%\[
%\max_{1\leq i < j\leq n}| d(x_i,x_j) - d(a_i,a_j) |,
%\]
%where the first $d$ is the metric symbol in the language and the second $d$ is the metric on $A$.

The theory of the Urysohn sphere is axiomatized by sentences which express an extension property characterizing the space we described in Example \ref{ury-introduced} 
We let the theory $T^{ext}$ consist of the sentences $\varphi_{A,B}$ defined by (recalling the notation from Example \ref{ex4.3})
\[
\sup_{\bar x} ( \inf_{y} D_B(\bar x,y) \dotminus D_A(\bar x)),
\]
where $n\geq 1$, $A = \{a_1,\ldots,a_n\}$ and $B = \{a_1,\ldots,a_n,b\}$ are $1$-bounded metric spaces, and $\bar x = x_1,\ldots,x_n$.

From our original characterization of the Urysohn sphere, it is clear that all the sentences $\varphi_{A,B}$ are satisfied in $U$.  

Before we show that $T^{ext}$ axiomatizes the theory of $U$, we first describe a different amalgamation construction than the one from Example \ref{ury-introduced}.  Suppose that $A$, $B = A \cup \{b\}$ and $C = A \cup \{c\}$ are all finite metric spaces with $b\not=c$ and with metrics $d_A, d_B$ and $d_C$ respectively for which $d_B$ and $d_C$ extend the metric $d_A$.  We wish to put a metric $d$ on $B \cup C$ which extends both $d_B$ and $d_C$; in other words, we must decide the value of $d(b,c)$. For the purposes of the argument to follow, we set
\[
d(b,c) = \max_{a \in A} |d_B(b,a) - d_C(c,a)|.
\]

\begin{exercise}
 Prove that $d$ is indeed a metric on $B\cup C$.   
\end{exercise}
We denote the resulting the resulting metric space $B \sqcup_A C$.

Now to see that $T^{ext}$ in fact axiomatizes the theory of $U$, it suffices to prove the stronger statement that any separable model of $T^{ext}$ is isomorphic to $U$.
So suppose that $\M$ is a separable model of $T^{ext}$.  One proves that $\M$ is isomorphic to $U$ using a standard back and forth argument.  We will concentrate on the inductive step and leave the rest of the details to the reader.  Consider a finite metric space $A=\{a_1,\ldots,a_n\} \subset \M$ and an abstract one point extension $B = A \cup \{b\}$ of $A$.  We wish to find an isomorphic copy of $B$ extending $A$ in $\M$.  First of all, we use the fact that $\varphi_{A,B}^\M=0$ to find some $b_0 \in \M$ such that $D^\M_B(a_1,\ldots,a_n,b_0) < 1/2$.   If we set $r_a = d(b,a)$ for $a \in A$, then we have that $|d(b_0,a) - r_a| < 1/2$ for all $a \in A$. 

Set $B_0:=\{a_1,\ldots,a_n,b_0\}$.  We now inductively build a sequence $(b_m)$ from $\M$ such that, setting $B_m :=\{a_1,\ldots,a_n,b_m\}$ and $\hat B_m:=B_m\sqcup_A B$, we have 

\[
D_{\hat{B}_m}(a_1,\ldots,a_n,b_m,b_{m+1}) < 1/2^{m+1}.
\]

% .Now we inductively build $B_m \subset \M$ with $B_m =\{a_1,\ldots,a_n,b_m\}$ such that
% \[
% D_{\hat{B}}(a_1,\ldots,a_n,b_m,b_{m+1}) < 1/2^n,
% \]
% where $\hat{B}$ is $B_m \sqcup_A B$ as described above. 

Note that we can indeed obtain $b_{m+1}\in \M$ by invoking the fact that
$\varphi_{B_m,\hat{B}_m}^\M=0$. This construction yields $|d(b_m,a) - r_a| < 1/2^{m+1}$ and also $d(b_m,b_{m+1}) < 3/2^{m+2}$.   It is now clear that $(b_m \colon m\in\bbN)$ is a Cauchy sequence in $\M$; letting $\tilde{b}\in \M$ be the limit of this sequence, we get that $d(\tilde{b},a) = r_a$ for all $a \in A$ and so the map $a_1,\ldots,a_n,b \mapsto a_1,\ldots,a_n,\tilde{b}$ is an isometry, yielding the desired copy of $B$ in $\M$.
\end{example}

We now turn to some examples that have more structure than just the metric symbol and are key to many of the articles in this volume.

\begin{example}[Hilbert space]\label{HS-axioms}
We now write sentences (or axioms) that are true in all (dissections of) Hilbert spaces using the language $\cL_{HS}$.  These sentences will in fact axiomatize the theory of (dissections of) Hilbert spaces.

The first set of axioms will express that we are dealing with a complex vector space.  Due to the presence of the sorts $S_n$, we will need infinitely many sentences to express the usual axioms.  For instance, for each $n\in \bbN$, we will have a sentence
\[
\sup_{x,y \in S_n} d_n(\lambda_n(x +_n y),\lambda_n(x) +_m \lambda_n(y)),
\]
where $\lambda \in \bbC$ and $S_m$ is the co-domain of $\lambda_n$.  Of course, if these sentences hold in our structure, then we will know that $\lambda(x + y) = \lambda x + \lambda y$ for all $x$ and $y$.  We leave it as an exercise to write out in a similar fashion all other sentences needed to guarantee that we have a complex vector space.

The second set of axioms expresses that we have an inner product and expresses its relationship with the metric on each sort.  For instance, we have the sentences
\[
\sup_{x \in S_n} | d_n(x,0)^2 - Re\langle x,x \rangle |.
\]
We leave it as an exercise to express the other axioms related to the inner product.

We need a third set of axioms which tell us that the inclusion maps between the sorts are isometric embeddings that preserve the algebraic operations.  First, we need to express that the inclusion maps themselves are compatible: for $n < m < k$, we include sentences
\[
\sup_{x \in S_n} d_k(i_{m,k}i_{n,m}(x), i_{n,k}(x)),
\]
which assert that $i_{m,k}i_{n,m} = i_{n,k}$  An instance of the inclusion maps preserving the algebraic structure are the sentences, for $n < m$ and $\lambda \in \bbC$,
\[
\sup_{x \in S_n} d_l(i_{k,l}(\lambda_n x) , \lambda_m i_{n,m}(x)),
\]
where $S_k$ is the co-domain of $\lambda_n$ and $S_l$ is the co-domain of $\lambda_m$.

Finally, we want to express that the sorts are interpreted correctly, that is, $S_n$ should be the ball of radius $n$.  We express this fact with two axioms:
\[
\sup_{x \in S_1} (d_1(x,0) \dotminus 1)
\]
and
\[
\sup_{x \in S_n} \min\{ 1 \dotminus \|x\|, \inf_{y \in S_1}d_n(y,i_{1,n}(y)\}.
\]
The first axiom says that everything in $S_1$ should have norm at most 1 and by scaling, everything in $S_n$ should have norm at most $n$.  The second axiom says that if something in $S_n$ has norm less than $1$, then it is in the closure of the image of the inclusion map from $S_1$.  Using the completeness of the sorts, we obtain that anything of norm at most 1 in $S_n$ arises from something in $S_1$.  Consequently, if these two axioms are satisfied (along with all the others above), then $S_1$ is the unit ball of a Hilbert space, and the direct limit of the sorts, via the embeddings, is a Hilbert space.  

If we denote by $T_{HS}$ the set of sentences listed above, then we conclude:
\begin{thm}
$T_{HS}$ axiomatizes the class of dissections of Hilbert spaces.  Moreover, models of $T_{HS}$ with homomorphisms as morphisms and the category of Hilbert spaces with homomorphisms as morphisms are equivalent via the dissection functor.
\end{thm}
\end{example}

\begin{example}[\cstar-algebras]\label{C*-axioms}
As in the previous example, we now write out sentences that axiomatize the class of (dissections of) \cstar-algebras in the language $\cL_{C^*}$ from Example \ref{C*ex}.

First of all we, need to express that the structure is a complex *-algebra.  This requires writing out a number of equations in the manner described above in the Hilbert space example.  We will not write out all the equations, but here is an example:
\[
\sup_{x,y \in S_n} d_n((xy)^*,y^*x^*).
\]
There is one such sentence for every $n \in \bbN$.  Of course, this set of sentences is expressing the relationship between multiplication and the adjoint in a \cstar-algebra, namely that $(xy)^* = y^*x^*$ for all $x$ and $y$.

We also want to say that the underlying structure is a normed algebra.  Again, we only write out a sample and leave it as an exercise to write out the rest.  It is tempting to write the triangle inequality as
\[
\forall_{x,y \in S_n} \|x + y\| \leq \|x\| + \|y\|
\]
but more formally it is written as
\[
\sup_{x,y \in S_n}\left ( \|x + y\| \dotminus (\|x\| + \|y\|) \right ),
\]
where $\|z\|$ is shorthand for $d(z,0)$. The latter sentence evaluates to $0$ in a $\cL_{C^*}$-structure precisely when the triangle inequality holds.

We also want our Banach *-algebra to be a \cstar-algebra.  The \cstar-identity is written as
\[
\sup_{x \in S_n} \left | \|x^*x\| - \|x\|^2 \right |,
\]
where we use the above convention for $\| \cdot \|$.  We include one such sentence for every $n \in \bbN$ although, in light of the other axioms, it would be enough to include this axiom only for $n = 1$.

As in the case of Hilbert spaces, there are now some sentences which must be included in order to make the formalism of continuous logic match with the normal operator algebraic setting.  To start with, we have the embedding maps $i_{m,n}$.  We can write out sentences which will express that these maps are isometric embeddings and they preserve the relevant algebraic structure (addition, multiplication and the adjoint).  We leave the formulation of these sentences to the reader.

We also want the sorts $S_n$ to be interpreted as the operator norm balls of radius $n$.  We need two sets of sentences.  First, we include sentences
\[
\sup_{x \in S_n} (\| x \| \dotminus n),
\]
which says that every element of the sort $S_n$ has norm at most $n$.  Finally, we record sentences  which will guarantee that all the elements of norm at most $n$ lie in $S_n$.  Recall from Example \ref{C*ex} the introduction of function symbols $\mathbf{q^n_k}$ with domain $S_n$ and co-domain $S_1$.  They are to represent polynomials $q^n_k$ and so we write sentences
\[
\sup_{x \in S_n} d_m(i_{1,m}(\mathbf{q^n_k}(x)), q^n_k(x)),
\]
where $S_m$ is the co-domain of the term $q^n_k$.  We will call the set of sentences that we have outlined in this example $T_{C^*}$.
By the argument after Lemma \ref{func-analytic}, we see that the $\cL_{C^*}$-structures which satisfy $T_{C^*}$ are precisely those of the form $D(A)$ for some \cstar-algebra $A$.  We conclude:
\begin{thm}
$T_{C^*}$ axiomatizes the class of dissections of \cstar-algebras.  Moreover the category of \cstar-algebras and the category of models of $T_{C^*}$ (with homomorphisms as morphisms) are equivalent via the dissection functor.
\end{thm}

\end{example}

\section{Ultraproducts}\label{sec-ultra}
\subsection{Ultraproducts of metric structures}
Fix a metric language $\cL = (\mathfrak{S, F, R})$, an index set $I$, an ultrafilter $\cU$ on $I$, and $\cL$-structures $\M_i$ for $i \in I$.  We construct the \emph{(metric) ultraproduct}\index{ultraproduct} $\N:=\prod_\cU \M_i$ as follows:
\begin{itemize}
\item For each sort $S \in \mathfrak{S}$, we let $S^\N$ be the metric space $\prod_\cU (S^{\M_i},d_i)$, where $d_i$ is the metric on $S^{\M_i}$.  We interpret the metric on $S^\N$ as $\lim_\cU d_i$. Recall that this is well-defined since all the metric spaces involved are uniformly bounded by the bound $B_S$.
\item For each $f \in \mathfrak{F}$ whose domain is $S_1 \times \ldots \times S_n$, we define $f^\N$ by
\[
f^\N([(a^1_i)]_\cU,\ldots,[(a^n_i)]_\cU) = [(f(a^1_i,\ldots,a^n_i))]_\cU,
\]
where $[(a^j_i)]_\cU \in S^\N_j$ for all $j = 1,\ldots,n$.  Again, since all interpretations of $f$ satisfy the same uniform continuity modulus, this definition of $f^\N$ is well-defined and $f^\N$ satisfies the uniform continuity modulus specified by the language.
\item For each $R \in \mathfrak{R}$, we define $R^\N$ by
\[
R^\N([(a^1_i)]_\cU,\ldots,[(a^n_i)]_\cU) = \lim_\cU R^{\M_i}(a^1_i,\ldots,a^n_i),
\]
where $[(a^j_i)]_\cU \in S^N_j$ for all $j = 1,\ldots,n$.  As in the case of function symbols, since the interpretation of $R$ is uniformly bounded in all $\M_i$ with the same uniform continuity modulus, this definition of $R^\N$ is well-defined and $R^\N$ satisfies the uniform continuity modulus specified by the language.
\end{itemize}
\begin{exercise}
Verify that $\prod_\cU \M_i$ is indeed an $\cL$-structure.
\end{exercise}

If all the $\M_i$ are the same $\cL$-structure $\M$, then we call the resulting ultraproduct the \emph{ultrapower}\index{ultrapower} and denote it by $\M^\cU$.

As in classical logic, it is important to be able to evaluate formulas in the ultraproduct.  Here is the appropriate version of \emph{\L o\' s' Theorem}\index{\L o\' s' Theorem}.
\begin{thm}
Suppose that $\M_i$ is an $\cL$-structure for all $i \in I$ and $\cU$ is an ultrafilter on $I$.  Set $\N = \prod_\cU \M_i$.  For an $\cL$-formula $\varphi(\bar x)$ and $\bar a = [(\bar a_i)]_\cU \in \N$ sorted the same as the variables $\bar x$, we have
\[
\varphi^\N(\bar a) = \lim_\cU \varphi^{\M_i}(\bar a_i).
\]
\end{thm}

\proof The proof is by induction on the formulation of basic formulas.  We only treat the case of $\varphi(\bar x) = \inf_y \psi(y,\bar x)$.

Suppose that $ \lim_\cU \varphi^{\M_i}(\bar a_i) > r$ for some $r \in \bbR$.  Set $X = \{ i \in I : \varphi^{M_i}(\bar a_i) > r \}$ and note that $X \in \cU$.  Fix $b = [(b_i)]_\cU \in \N$. For each $i \in X$, we have $\psi^{\M_i}(b_i,\bar a_i) > r$.  By induction, we have $\psi^\N(b,\bar a) \geq r$. Since this was true for all $b \in \N$, we have $\varphi^\N(\bar a) \geq r$ and thus $\varphi^\N(\bar a) \geq \lim_\cU \varphi^{\M_i}(\bar a_i)$.

Now suppose that $s  = \varphi^\N(\bar a) > r$ but $\lim_\cU  \varphi^{\M_i}(\bar a_i) \leq r$.  Pick $\e > 0$ small enough so that $r + \e < s$ and let $X = \{ i \in I \colon \varphi^{\M_i}(\bar a_i) < r + \e \} \in \cU$.  For $i \in X$, choose $b_i \in \M_i$ so that $\psi^{\M_i}(b_i,\bar a_i) < r + \e$.  Set $b = [(b_i)]_\cU$.  By induction we have $\psi^\N(b,\bar a) \leq r+ \e < s$, which is a contradiction. \qed

\begin{exercise}
Verify the remaining cases in the proof of \L os's theorem including generalizing to arbitrary formulas from basic formulas.
\end{exercise}

\subsection{Applications of the ultraproduct}
We say that a set of sentences $\Gamma$ is \emph{finitely satisfiable}\index{finitely satisfiable} if every finite subset of $\Gamma$ is satisfiable.  We say that $\Gamma$ is \emph{approximately finitely satisfiable}\index{approximately finitely satisfiable} if, for every finite subset $\Gamma_0 \subset \Gamma$ and every $\e > 0$, there is an $\cL$-structure $\M$ so that $|\varphi^\M| \leq \e$ for every $\varphi \in \Gamma_0$.  Equivalently, $\Gamma$ is approximately finitely satisfiable if the set of sentences $\{ |\varphi| \dotminus \e \colon \varphi \in \Gamma,\  \e>0 \}$ is finitely satisfiable.

\begin{thm}[Compactness Theorem]
For a metric language $\cL$ and a set of $\cL$-sentences $\Gamma$, the following are equivalent:
\begin{enumerate}
\item $\Gamma$ is satisfiable.
\item $\Gamma$ is finitely satisfiable.
\item $\Gamma$ is approximately finitely satisfiable.
\end{enumerate}
\end{thm}

\proof Clearly the conditions decrease in strength.  To show that (3) implies (1), consider $I$ to be the set of finite subsets of $\Gamma \cup \bbN$ (disjoint union) which have non-empty intersection with both $\Gamma$ and $\bbN$.  If $i = F \cup X$ is an element of $I$ with $F \subset \Gamma$ and $X \subset \bbN$, let $M_i$ be an $\cL$-structure such that $|\varphi^{M_i}| \leq 1/m$ for all $\varphi \in F$, where $m = \max\{ n \colon n \in X\}$.  For $\varphi \in \Gamma$ and $n \in \bbN$, let $O_{\varphi,n}$ be the set of $i \in I$ which contain both $\varphi$ and $n$.  Notice that $Y = \{ O_{\varphi,n} \colon \varphi \in \Gamma, n \in \bbN\}$ has the finite intersection property. Let $\cU$ be an ultrafilter on $I$ which contains $Y$.  Consider the ultraproduct $\N = \prod_\cU M_i$.  For every $\varphi \in \Gamma$ and every $n \in \bbN$, we have $O_{\varphi,n} \in \cU$, so by the \L o\' s Theorem, we have $|\varphi^\N|\leq \frac{1}{n}$ for all $n$ and hence $\varphi^\N=0$.  It follows that $\N$ satisfies $\Gamma$. \qed

We can use the compactness theorem and the method of expansion by constants to give a version of the upward L\" owenheim-Skolem theorem in the continuous setting.  We say that a metric structure $\M$ is \emph{compact} if $S^\M$ is a compact metric space for each sort $S$ in the language.
\begin{thm}
Suppose that $\M$ is a non-compact $\cL$-structure and $\chi(Th(\M),\cL) \leq \lambda$.  Then $\M$ has an elementary extension $\N$ such that $\chi(\N) = \lambda$.
\end{thm}

\proof Since $\M$ is not compact, there is a sort $S$ and an $\e > 0$ such that the collection of $\e$-balls of $S^\M$ does not have a finite subcover.  Now consider the language $\cL'$ obtained from $\cL$ by introducing new constant symbols $c_\alpha$ of sort $S$ for $\alpha < \lambda$.  Consider the $\cL'$-theory $T'$ made up of $Elem(\M)$ together with
\[
\{ d(c_\alpha,c_\beta) \geq \e : \alpha < \beta < \lambda \}.
\]
Note that $T'$ is finitely satisfiable and $\chi(T',\cL') = \lambda$.  By the compactness theorem and the downward L\" owenheim-Skolem theorem, $T'$ has a model with density character $\lambda$.  Taking the reduct of this model to $\cL$, we get the required model of $T$ of density character $\lambda$. \qed

\subsection{Elementary equivalence and ultraproducts}
We begin this section with an exercise.
\begin{exercise}
For $\cL$-structures $\M$ and $\N$, we have that $\M \equiv \N$ if and only if there is an $\cL$-structure $\cal P$ such that $\M$ and $\N$ both elementarily embed into $\cal P$.
\end{exercise}

%This fact can be significantly tightened up.
This exercise can be significantly strengthened:
\begin{thm}[Keisler-Shelah]
For $\cL$-structures $\M$ and $\N$, we have that $\M \equiv \N$ if and only if there are ultrafilters $\cU$ and $\cV$ such that $\M^\cU \cong \N^\cV$.
\end{thm}
The proof in the classical case can be found in \cite{Keisler2, Sh}.  In the continuous case, a proof appears in \cite{GK} while a proof in the positive bounded case appears in \cite{HenIov}.
It is important to notice that in the Keisler-Shelah theorem, the index sets for the ultrafilters need not be countable.  

In the continuous setting, it is often difficult to explicitly write out a sentence that distinguishes the theories of two structures.  The approach of showing that no two ultrapowers are isomorphic has been used to good effect in separating theories of II$_1$ factors (see \cite{BCI}, \cite{CIK} and Section 4 in the author's article with Goldbring in this volume).

EF-games were discussed in the previous section. There is also an infinite version of an EF-game involving ultrapowers which is sometimes easier to work with.  We describe this game only in the case where $\cL$ is separable. Given two $\cL$-structures $\M$ and $\N$, the $\infty$-$EF(\M,\N)$-game is a two player game between Player I and Player II which now goes on for $\omega$ many rounds.  The players play as in a finite EF-game, this time creating two sequences $(a_i : i \in \bbN)$ in $\M$ and $(b_i : i \in \bbN)$ in $\N$.  Player II wins the game if the map $a_i \mapsto b_i$, for $i \in \bbN$, extends to an isomorphism of $\cL$-substructures of $\M$ and $\N$ respectively.

\begin{thm}
Assume $\cL$ is separable. Then for two $\cL$-structures $\M$ and $\N$, the following are equivalent:
\begin{enumerate}
\item $\M \equiv \N$.
\item For some (equiv. any) non-principal ultrafilter $\cU$ on $\bb N$, player II has a winning strategy for the $\infty$-$EF(\M^\cU,\N^\cU)$-game.
\end{enumerate}
\end{thm}

\begin{exercise}
Prove the previous theorem.
\end{exercise}
\subsection{Elementary classes}

  We now give a characterization of elementary classes which parallels a result from classical first order logic. Fix a class $\cal C$ of $\cL$-structures.    We say that $\cal C$ is \emph{closed under ultraroots}\index{ultraroots} if for any ultrafilter $\cU$, whenever $\N^\cU \in \cal C$, then $\N \in \cal C$.
\begin{thm}\label{elementaryclasstheorem} Suppose that $\cal C$ is a class of $\cL$-structures.  The following are equivalent:
\begin{enumerate}
\item $\cal C$ is an elementary class.
\item $\cal C$ is closed under isomorphisms, ultraproducts, and elementary submodels.
\item $\cal C$ is closed under isomorphisms, ultraproducts, and ultraroots.
\end{enumerate}
\end{thm}

\proof Clearly (1) implies both (2) and (3).   We  first prove (2) implies (1) and then modify the proof to get (3) implies (1).  Thus, assume (2) and let $T$ be the theory of $\cal C$, that is, let $T$ be the set of all sentences $\varphi$ such that $\varphi^\N = 0$ for all $\N \in \cal C$.  We show that $\cal C$ is the collection of models of $T$.  It is clear that all elements of $\cal C$ model $T$.  On the other hand, fix a model $\M$ of $T$. By (2), it suffices to elementarily embed $\M$ into an ultraproduct of structures from $\cal C$.  To this end, suppose $\varphi(c_{\bar m})$ is in $Elem(\M)$ for some sequence of constants $c_{\bar m}$ with $\bar m \in \M$.  Now consider $\inf_{\bar x} |\varphi(\bar x)|$, an $\cL$-sentence whose value in $\M$ is 0.  Fix $\epsilon>0$.  Then there is $\N_{\varphi,\e}\in \cal C$ which satisfies $\inf_{\bar x}| \varphi(\bar x)| < \e$ (for otherwise the sentence $\epsilon \dotminus \inf_{\bar x}|\varphi(\bar x)|$ belongs to $T$, contradicting the fact that $\M\models T$). By taking an ultraproduct over all $\N_{\varphi,\e}$ as $\varphi$ and $\e$ vary, we can elementarily embed $\M$ into an ultraproduct of elements of $\cal C$, as desired.  

To obtain (3) implies (1), notice that the previous argument showed (assuming that $\cal C$ is merely closed under isomorphisms and ultraproducts) that $\M \equiv \N$ for some $\N$ in $\cal C$ (namely an ultraproduct of structures from $\cal C$).  By the Keisler-Shelah theorem, $\M^\cU \cong \N^\cV$ for some ultrafilters $\cU$ and $\cV$. Assuming (3), we have that $\M \in \cal C$. \qed

Although the semantic characterization of elementary class is the same in continuous logic as in classical logic, it plays a somewhat bigger role in the continuous setting as it is often difficult to explicitly write out sentences that axiomatize a given class of metric structures. For instance, although the class of W$^*$-probability spaces (see \cite{Ando} for the definition) was axiomatized in a complicated language in \cite{D}, it was relatively easy to show using the above semantic characterization that this class was elementary in a more natural language \cite{GHS}; an explicit axiomatization in this language will only finally appear in \cite{AGH}.

\begin{exercise}
Referring back to Example \ref{C*ex} of the class of \cstar-algebras and their dissections as metric structures, we want to see that the ultraproduct in the metric structure sense and the \cstar-algebra sense match up.  Recall that 
if $(A_i:i\in I)$ is a family of \cstar-algebras and $\cU$ is an ultrafilter on $I$, then
we set $\prod_I^\infty A_i$ to be
\[
\{ \bar a \in \prod_I A_i : \text{there is } B>0 \text{ so that } \|a_i \| \leq B \text{ for all } i \in I\}
\]
and
\[
J = \{ \bar a \in \prod_I^\infty A_i: \lim_{\cU} \| a_i \| = 0 \}.
\]
The C*-algebraic ultraproduct of $(A_i)$ with respect to $\cU$ is then defined to be 
\[
\prod_\cU A_i = \prod_I^\infty A_i / J.
\]
Prove that $D(\prod_\cU A_i) = \prod_\cU D(A_i)$, where the ultraproduct appearing on the right-hand side is the ultraproduct of metric structures.
\end{exercise}

This exercise, together with Example \ref{C*ex}, shows, by Theorem \ref{elementaryclasstheorem}, that the class of dissections of \cstar-algebras is an elementary class, categorically equivalent to the class of \cstar-algebras via the dissection functor which even preserves ultraproducts.  Of course, we already know that the class of dissections is an elementary class because they are exactly the $\cL_{C^*}$-structures  which satisfy the theory $T_{C^*}$ from Example \ref{C*-axioms}.

\section{Types}\label{sec-types}

\subsection{The type space}
Suppose that $T$ is a theory in some language $\cL$.  Fix a sequence of sorted variables $\bar x$.  The \emph{space of $\bar x$-types of $T$}\index{type space}, denoted $S_{\bar x}(T)$, is the set of functions of the form $tp^\M(\bar a)(\varphi(\bar x)) := \varphi^\M(\bar a)$, where $\M$ is a model of $T$, $\bar a \in \M$ (sorted in the same way as the variables in $\bar x$), and $\varphi$ is an $\cL$-formula whose free variables are among $\bar x$.  $tp^\M(\bar a)$ is called the \emph{complete type of $a$ in $\M$}\index{complete type} (complete meaning that the value of every formula $\varphi$ is completely specified).  We use letters such as $p$ and $q$ to denote complete types.  If $p\in S_{\bar x}(T)$ and $\bar a\in \M$ is such that $tp^\M(\bar a)=p$, then we say that $\bar a$ \emph{realizes} $p$ and that $p$ is \emph{realized in $\M$}.  Just as with complete theories, complete types are determined by knowing the kernel of the type functional:
\[
tp^\M(\bar a)(\varphi(\bar x)) = r \text{ iff } \varphi(\bar x) - r \in ker(tp^\M(\bar a)).
\]
It is convenient to talk about partial types as well.  Suppose $\Sigma$ is a set of $\cL$-formulas in free variables $\bar x$ and $p$ is a complete type in the variables $\bar x$.  Then the corresponding \emph{partial ($\Sigma$)-type}\index{partial type} is the function $p\restriction \Sigma$.
If $\Sigma$ is the set of all quantifier-free formulas in the free variables $\bar x$, we call the $\Sigma$-type a \emph{quantifier-free type}.

As in the classical case, the type space $S_{\bar x}(T)$ carries a natural topology.  For $\varphi \in Form_{\bar x}(T,\cL)$ and real numbers $r < s$, set 
\[
O^T(\varphi,(r,s)) = \{ p \in S_{\bar x}(T) \colon p(\varphi)\in (r,s) \}.
\]
  
  \begin{exercise}

  \
  
      \begin{enumerate}
          \item Prove that the set of all $O^T(\varphi, (r,s))$ is a base for a Hausdorff topology on $S_{\bar x}(T)$, called the \emph{logic topology}\index{logic topology}.
          \item Prove that the logic topology on $S_{\bar x}(T)$ is compact.  (Hint:  use the compactness theorem.)
          \item Prove that if $\pi(\bar x)$ is a partial type, then the set $\{p\in S_{\bar x}(T) \ : \ p \text{ extends }\pi\}$ is a closed subset of $S_{\bar x}(T)$.
          \item In fact, show that every non-empty closed in the logic topology on $S_{\bar x}(T)$ has the form $\{p\in S_{\bar x}(T) \ : \ p \text{ extends }\pi\}$, where $\pi$ is a partial type.  (Hint: Show that those $p\in S_{\bar x}(T)$ for which $p(\varphi) \leq r$ are the same as those which extend the partial type $\pi$ where $\pi(\varphi \dotminus r) = 0$.)
          \end{enumerate}
      \end{exercise}

If the theory $T$ is complete, the type space also has a natural \emph{metric topology}\index{metric topology}: for $p,q \in S_{\bar x}(T)$, we set $d(p,q)$ to be the infinum of $d(\bar a,\bar b)$, where $p = tp^\M({\bar a})$ and $ q = tp^\M({\bar b})$ for some model $\M$ of $T$ and some $\bar a, \bar b \in \M$ . The completeness of the theory guarantees that $d(p,q)$ is defined (not infinite) for every pair of types $p$ and $q$.  The symmetry and reflexivity of $d(p,q)$ are clear.  To see that the triangle inequality holds, fix $p,q,r \in S_{\bar x}(T)$. By compactness, we can find models $\M$ and $\N$ of $T$ and $\bar a, \bar b \in \M$, $\bar b', \bar c \in \N$ with $p=tp^\M(\bar a)$, $q=tp^\M(\bar b)=tp^\N(\bar b')$, $r=tp^\N(\bar c)$, $d(p,q) = d(\bar a,\bar b)$ and $d(q,r) = d(\bar b',\bar c)$.  We leave it as an exercise that by taking a common elementary extension of $\M$ and $\N$, we can assume that $\M = \N$ and $\bar b = \bar b'$.  By the triangle inequality (in $\M$), $d(p,q) + d(q,r) \geq d(\bar a,\bar c) \geq d(p,r)$.

\begin{exercise}
\begin{enumerate}
    \item Show that the metric topology refines the logic topology on $S_{\bar x}(T)$.
    \item Show that $S_{\bar x}(T)$ is complete with respect to the metric topology.
\end{enumerate}    
\end{exercise}

We will see later that there are model theoretic consequences when these two topologies agree (even locally).

\subsection{Some uses of the Stone-Weierstrass theorem in continuous logic}\label{SW}

To understand the extension from basic formulas to more general formulas in a more abstract fashion, we remind the reader of the background necessary for the Stone-Weierstrass theorem.

Suppose that $X$ is a topological space and set $C(X,\bbR)$ to be the set of continuous functions from $X$ to $\bbR$.  $C(X,\bbR)$ has a natural real algebra structure by considering pointwise addition and multiplication together with scalar multiplication by elements of $\bbR$.  It is also naturally topologized by uniform convergence, that is, for $f \in C(X,\bbR)$ and $\e > 0$, the set
\[
\{ g \in C(X,\bbR) \colon |f(x) - g(x)| < \e \text{ for all } x \in X \}
\]
is a basic open set. When $X$ is compact, $C(X,\bbR)$ is metrizable by the $\sup$-norm. For $A \subseteq C(X,\bbR)$, we say that $A$ \emph{separates points} if for all $x,y \in X$, there is $f \in A$ such that $f(x) \neq f(y)$.

\begin{thm}[Stone-Weierstrass theorem]
Suppose that $X$ is a compact, Hausdorf space and $A$ is a real subalgebra of $C(X,\bbR)$.  If the constant function $1$ belongs to $A$, then $A$ is dense in $C(X,\bbR)$ if and only if $A$ separates points.
\end{thm}

If $\varphi(\bar x)$ is an $\cL$-formula and $T$ is an $\cL$-theory, then there is a natural evaluation map $f_\varphi \colon S_{\bar x}(T) \rightarrow \bbR$ given by $f_\varphi(p) = p(\varphi)$. With the logic topology on the type space $S_{\bar x}(T)$, it is easy to see that $f_\varphi$ is continuous. We want to use the Stone-Weierstrass theorem to give another understanding of formulas in continuous model theory.  In fact, this alternative perspective is really the reason behind extending the class of formulas.

\begin{thm}\label{formula-thm}
Suppose that $T$ is a theory, $\bar x$ a sequence of variables, and $X = S_{\bar x}(T)$.  Then the evaluation map $ev$ from $Form_{\bar x}(T,\cL)$ to $C(X,R)$ given by $ev(\varphi) = f_\varphi$ is an isometry which commutes with all continuous real-valued functions, that is, if $f:\bbR^n \rightarrow \bbR$ is a continuous function and $\varphi_1,\ldots,\varphi_n \in Form_{\bar x}(T,\cL)$, then $$f\circ(ev(\varphi_1),\ldots,ev(\varphi_n)) = ev(f(\varphi_1,\ldots,\varphi_n)).$$
\end{thm}

\proof Note that the range of the evaluation map restricted to basic formulas is a real subalgebra containing 1 and which separates points.  It follows then from the  Stone-Weierstrass theorem that the evaluation map is onto.  The map is norm-preserving since if $\varphi(\bar x)$ is a basic formula such that $\|\varphi\|_T = r$, then on the one hand, $\|ev(\varphi)\| \leq r$ in $C(X,\bbR)$, while by compactness, we can find some $p \in S_{\bar x}(T)$ such that $ev(|\varphi|)(p) = r$, hence $\|ev(\varphi)\| = r$.  Commutation with continuous functions is immediate. \qed

We need a slight strengthening of the previous theorem for real algebra systems (see Definition \ref{real-alg}).  The proof follows from compactness and the Stone-Weierstrass theorem.

\begin{prop}\label{gen-SW}
Suppose that $T$ is an $\cL$-theory and $A$ is real algebra system of $T$-formulas.  Let $X$ be the set of restrictions of elements in $S_{\bar x}(T)$ to formulas in $A_{\bar x}$.  Then
\begin{enumerate}
\item $X$ is compact, and
\item if the range of the evaluation map from $A_{\bar x}$ to $C(X,\bbR)$ separates points then the evaluation map is an isometry that commutes with all continuous functions.
\end{enumerate}
\end{prop}

\begin{exercise}
    Prove the previous proposition.
\end{exercise}

As mentioned before, three places where the previous proposition will be useful are for the following real algebra systems:
\begin{enumerate}
\item $A_{\bar x}$ is the closure of the basic $\cL'$-formulas in $Form_{\bar x}(T,\cL)$ for $\cL'$, a sublanguage of $\cL$, and $T$, an $\cL$-theory.
\item  $A_{\bar x}$ the set of quantifier-free $\cL$-formulas in the free-variables $\bar x$.
\item $A_{\bar x}$ the set of prenex formulas in $\cL$ in the free variables $\bar x$.  Note that prenex formulas in $\cL$ are dense in the set of all $\cL$-formulas. (See Section 6 of \cite{BBHU} for details.)
\end{enumerate}

\subsection{Saturation}

Suppose that $\M$ is an $\cL$-structure and $A \subset \M$.  We say that $p$ is a \emph{(partial) type over $A$}\index{partial type} if $p$ is a (partial) type in $\cL_A$ with respect to the theory of $\M_A$.  Consequently, if $p$ is a partial type over $A$, there is an elementary extension $\N$ of $\M$ and $\bar b \in \N$ realizing $p$, that is, if $\varphi(\bar x,\bar a)$ is an $\cL_A$-formula (in the domain of $p$), then $p(\varphi(\bar x,\bar a)) = \varphi^\N(\bar b,\bar a)$.

\begin{defn}
For $\kappa$ an infinite cardinal, we say that an $\cL$-structure $\M$ is \emph{$\kappa$-saturated}\index{$\kappa$-saturated} if, whenever $A \subseteq \M$ is such that $\chi(A) < \kappa$ and $p$ is a type over $A$, then $p$ is realized in $\M$.  $\M$ is said to be \emph{saturated}\index{saturated} if it is $\chi(\M)$-saturated.
\end{defn}

As in the case of classical first order logic, saturated models typically only exist under additional set theoretic assumptions.  For instance, if $2^\kappa = \kappa^+$, then the following proposition yields that if $\chi(\M) = \kappa$, then $\M$ has a saturated elementary extension of density character $\kappa^+$.

\begin{prop}
Suppose that $\M$ is an $\cL$-structure and $\chi(\cL) \leq \chi(\M) = \kappa$.  Then there is an elementary extension $\N$ of $\M$ such that $\N$ is $\kappa^+$-saturated and $\chi(\N)\leq 2^\kappa$.
\end{prop}

\proof We construct $\N$ as a union of an increasing chain $(\M_\alpha)_{\alpha < \kappa^+}$ such that:
\begin{enumerate}
\item $M_0 = \M$,
\item if $\delta < \kappa^+$ is a limit ordinal, then $\M_\delta$ is the completion of $\bigcup_{\alpha < \delta} \M_\alpha$ (that is, we take the direct limit of the $\M_\alpha$'s as described in Proposition \ref{direct-limit}),
\item $\chi(\M_\alpha) \leq 2^\kappa$ for all $\alpha < \kappa^+$, and 
\item for every $\alpha < \kappa^+$, $\M_\alpha \prec M_{\alpha + 1}$ and for every $A \subseteq \M_\alpha$ with $\chi(A) \leq \kappa$, every type $p$ over $A$ is realized in $\M_{\alpha + 1}$.
\end{enumerate}

Only the last point needs comment.  There are $2^\kappa$ many subsets of size $\kappa$ of a set of size $2^\kappa$ and at most $2^\kappa$ many types over a set of density character $\kappa$.  Consequently, we need to realize $2^\kappa$ many types, which one can do one at a time with a chain of length $2^\kappa$ since any individual type can be realized in an elementary extension of density character at most $2^\kappa$. \qed

One model-theoretic reason for working with ultraproducts is that they are always somewhat saturated:

\begin{exercise}
Suppose that $\cL$ is a separable language, $\M_n$ is an $\cL$-structure for all $n \in \bbN$, and $\cU$ is a non-principal ultrafilter on $\bbN$.  Show that $\prod_\cU \M_n$ is $\aleph_1$-saturated.
\end{exercise}

%An ultrafilter $\cU$ on $\kappa$ is called regular if there there are $X_\alpha \in \cU$ for $\alpha < \kappa$ such that for all $i \in I$, $\{ \alpha < \kappa : i \in X_\alpha\}$ is finite.  For every $\kappa$, there always exists a regular ultrafilter on $\kappa$.  
Ultraproducts can also produce structures of higher saturation as the following proposition indicates.  The ultrafilters used in the proof of this proposition are called \emph{good ultrafilters} (see \cite{G-ultra}). A proof of the existence of good ultrafilters can be found in \cite{Kunen} and the saturation of the ultraproduct can be found in \cite{Keisler2}; see also \cite[Chapter 8]{G-ultra}.

\begin{prop}
For every infinite cardinal $\kappa$, there is an ultrafilter $\cU$ on $\kappa$ such that if
$\cL$ has density character at most $\kappa$ and $\M_\alpha$ is an $\cL$-structure for all $\alpha < \kappa$, then $\prod_\cU \M_\alpha$ is $\kappa^+$-saturated.
\end{prop}

\subsection{Beth definability}

We now prove the Beth Definability theorem in continuous logic. This theorem gives us an abstract criterion for recognizing relations that can be expressed by a formula. As we will see, this theorem plays a wider role in continuous logic than in classical logic due to the difficulty in identifying which relations are exactly formulas in the continuous setting.
\begin{thm}[Beth Definability]\label{Beth}
Suppose that $\cL' \subseteq \cL$ are two continuous languages with the same sorts.  Further, suppose $T$ is an $\cL$-theory.  If the forgetful functor
$F:Mod(T) \rightarrow Str(\cL')$ given by $F(\M) = \M \restriction \cL'$ is an equivalence of categories onto the image of $F$, then every $\cL$-formula is $T$-equivalent to an $\cL'$-formula.
\end{thm}

\proof Suppose that $\varphi(\bar x)$ is an $\cL$-formula.  For $\e > 0$, consider the set of sentences $\Sigma_\e$ in the language $\cL$ with two tuples of constants $\bar a$ and $\bar b$ added consisting of $T$ together with 
\[
\{ | \psi(\bar a) - \psi(\bar b) | \dotminus \frac{1}{n} : \psi(\bar x) \text{ is an $\cL'$-formula, } n \in \bbN \} \cup \{ \e \dotminus |\varphi(\bar a) - \varphi(\bar b)|\}.
\]
If $\Sigma_\e$ is unsatisfiable for all $\e$, then by compactness, $\varphi$ defines a continuous real-valued function on the $\cL'$-restrictions of $\cL$-types in the variables $\bar x$ and so by Proposition \ref{gen-SW} applied to the real algebra of $\cL'$-formulas, $\varphi$ is $T$-equivalent to an $\cL'$-formula, which is the desired conclusion.

Based on the previous paragraph, we now suppose, towards a contradiction, that $\Sigma_\e$ is satisfiable for some $\e > 0$.  We can thus find a model $\M$ of $T$ and $\bar a, \bar b \in \M$ such that $\bar a$ and $\bar b$ have the same $\cL'$-type but $\varphi^\M(\bar a) \neq \varphi^\M(\bar b)$.  

Take an $|M|^+$-saturated elementary extension $\N$ of $\M$. 
%We may think of $\M$ as an $\cL$-substructure of $\N$ via the diagonal embedding.  
Let $\M'$ and $\N'$ be the reducts of $\M$ and $
N$ to $\cL'$.  Since $\N'$ is also $|\M'|^+$-saturated, we can find an $\cL'$-elementary map $\rho\colon\M' \rightarrow N'$ such that $\rho(\bar b) = \bar a$.  Since $F$ is an equivalence of categories, $\rho$ is also an $\cL$-elementary map, which contradicts that $\varphi^\M(\bar a) \neq \varphi^\M(\bar b)$. \qed

\section{Definable Sets}\label{sec-def}

\subsection{Introducing definable sets}
We will now look at one of the biggest differences between continuous model theory and classical model theory.  In the classical setting, it is often advantageous to be able to quantify over a subset of the universe which is described by a formula.  For instance, if $\varphi(x)$ is a formula expressing some property, it is straightforward to also express "for all $x$, if $\varphi(x)$ then ..." or "there exists $x$ such that $\varphi(x)$ and ...".  In the continuous setting, where formulas take on real values, the subsets over which one can quantify is a delicate matter.  

Suppose that $\M$ is an $\cL$-structure and $\varphi(\bar x)$ is an $\cL$-formula.  The \emph{zero-set of $\varphi$ in $\M$}\index{zero-set}, written $Z^\M(\varphi)$, is the set $\{ \bar a \in \M \colon \varphi^\M(\bar a) = 0 \}$.  Naively one might think that one can quantify over zero sets of formulas, but as we will soon see, that is not quite correct.

Suppose that $T$ is an $\cL$-theory and $\bar S$ is a tuple of sorts from $\cL$.  A \emph{$T$-functor (for $\bar S$)}\index{$T$-functor} is a functor $X \colon Mod(T) \rightarrow Met$, where $Met$ is the category of metric spaces with isometric embeddings as morphisms, $X(\M)$ is a closed subspace of $\bar S(\M)$, and if $\rho:\M \rightarrow \N$ is an elementary map in $Mod(T)$, then $X(\rho)$ is the restriction of $\rho$ to $X(\M)$.

\begin{example}
A natural example of a $T$-functor is $X_\varphi$ for some formula $\varphi$, where $X_\varphi(\M) = Z(\varphi^\M) = \{\bar a \in \M \colon \varphi^\M(\bar a) = 0\}$, the zero set of $\varphi$.
\end{example}

\begin{thm}\label{def-thm} Suppose that $T$ is an $\cL$-theory and $X$ is a $T$-functor for $\bar S$. The following are equivalent:
\begin{enumerate}
%\item $X$ is definable.
\item There is a $T$-formula $\varphi(\bar x)$ such that for any model $\M$ of $T$ and any $\bar a\in \M$, 
\[
\varphi^\M(\bar a) = d(\bar a,X(\M)).
\]
\item For any $T$-formula $\psi(\bar x,\bar y)$, there is a $T$-formula $\chi(\bar y)$ such that for any model $\M$ of $T$ and any $\bar a \in \M$,
\[
\chi^\M(\bar a) = \inf \{ \psi^\M(\bar b,\bar a) \colon \bar b \in X(\M) \}
\] (or equivalently the analogous statement replacing $\inf$ by $\sup$).

\item For every $\e > 0$, there is a $T$-formula $\varphi(\bar x)$ and $\delta >0$ such that, for all models $\M$ of $T$:
\begin{itemize}
    \item  $X(M)\subseteq Z(\varphi^\M)$, and 
    \item For all $\bar a \in \M$, if $\varphi^\M(\bar a)< \delta$, then $d(\bar a,X(\M)) \leq \e$.
\end{itemize}

\item For every $\e > 0$, there is a basic $\cL$-formula $\varphi(\bar x)$ and $\delta > 0$ such that, for all models $\M$ of $T$:
\begin{itemize}
    \item $X(\M) \subseteq Z(\varphi^\M)$, and
    \item For all $\bar a \in \M$, if $\varphi^\M(\bar a) < \delta$, then $d(\bar a,X(\M)) \leq \e$.
\end{itemize}

\item For any set $I$, family of $\cL$-structures $\M_i$ for $i \in I$, and ultrafilter $\cU$ on $I$, 
\[
X(\N) = \prod_\cU X(\M_i), \text{ where }\N = \prod_\cU \M_i.
\]
\end{enumerate}
\end{thm}

\begin{comment}\label{definablecomment} 
Before we begin the proof of Theorem \ref{def-thm}, we want to address one small point.  What do we do with the empty functor, that is, the functor $X$ which satisfies $X(\M) = \emptyset$ for all models $\M$ of $T$?  The empty functor certainly satisfies condition (5) above; what about the other four conditions?  In condition (1) we must decide what the distance is to the empty set - it should be the maximum possible value of the metric (and hence is given by a formula).  In (2), the $\inf$ is taken over an empty set hence equals the maximum possible value of $\psi$ - again, a formula.  In (3) and (4), we are measuring the distance to the empty set and so we again take this to be the upper bound on the metric.  Consequently, the empty functor satisfies the equivalent conditions of Theorem \ref{def-thm} when conditions (1)-(4) are suitably interpreted.  Moreover, if $X$ satisfies any of the conditions of the theorem and $X(\M)$ is empty for some model $\M$ of $T$, then it is empty for all models $\M$ of $T$ (and thus satisfies all of the other conditions in the theorem).
\end{comment}

\proof[Proof of Theorem \ref{def-thm}] 
Based on the previous comment, we can assume, without loss of generality, that $X(\M)\neq\emptyset$ for all models $\M$ of $T$.

The proof of (1) implies (2) can be found in \cite[Theorem 9.17]{BBHU}.  We only remark that if $d_X(x)$ is a formula which expresses the distance from $x$ to $X$ in any model of $T$ and $\alpha$ is a uniform continuity modulus for the formula $\psi$, then in any model $\M$ of $T$ (remembering that $X(\M)$ is not empty), $ \inf \{ \psi^\M(\bar b,\bar a) \colon \bar b \in X(\M) \}$ equals
\[
\inf_{\bar b\in \M} (\psi^\M(\bar b,\bar a) + \alpha(d^\M_X(\bar b))).
\]

It is easy to see that (2) implies (3).  For (3) implies (4), it is enough to remember that $T$-formulas are uniformly approximated by basic $\cL$-formulas. To see that (4) implies (5), we first show that $\prod_\cU X(\M_i)\subseteq X(\N)$.  To this end, fix $\e>0$ and let $\varphi$ and $\delta$ be as in (4).  If $[(a_i)]_\cU\in \prod_\cU X(\M_i)$, then $\varphi^{\M_i}(a_i)=0$ for all $i\in I$, hence $\varphi^\N([a_i]_\cU)=0$ and thus $d([(a_i)]_\cU,X(\N))\leq \epsilon$.  Since $\e>0$ was arbitrary, we have that $[(a_i)]_\cU\in X(\N)$.  Now suppose, towards a contradiction, that there is $[(a_i)]_\cU\in X(\cal N)\setminus \prod_\cU X(\cal M_i)$.  Since $\prod_\cU X(\cal M_i)$ is closed in $\cal N$, there is $\epsilon>0$ such that $d([(a_i)]_\cU,\prod_\cU X(\cal M_i))>\epsilon$.  Let $\varphi$ and $\delta>0$ be as in (4).  Since $X(\cal N)\subseteq Z(\varphi^\cal N)$, we have that $\varphi^{\M_i}(a_i)<\delta$ for $\cU$-almost all $i\in I$, hence $d(a_i,X(\cal M_i))\leq \epsilon$ for $\cU$-almost all $i\in I$, which is a contradiction.

It remains to show that (5) implies (1).  Let $\cal{C}$ be the class of structures 
\[
\{ (\M,d_X^\M) \colon \M \models T \},
\]
where $d_X^\M$ measures the distance from an element to $X(\M)$.  Note that measuring the distance to a closed set is indeed an allowable relation in continuous logic: the bound is the maximum of the bounds on the metrics in $\bar S$ and the modulus of uniform continuity is the identity function.  By Theorem \ref{elementaryclasstheorem}, condition (5) guarantees that $\cal C$ is an elementary class.  The forgetful functor from $\cal C$ onto $Mod(T)$ is an equivalence of categories since any model of $T$ extends uniquely to an object in $\cal C$.  By the Beth Definability theorem, we conclude that $d_X$ is in fact equivalent to a $T$-formula.   \qed

If $X$ is a $T$-functor (for $\bar S)$ satisfying any of the equivalent conditions of the previous theorem, then we say that $X$ is a \emph{$T$-definable set}\index{$T$-definable}\index{definable set}. By condition (2) of the previous theorem, $(T$-)definable sets are exactly those sets in metric structures over which one can quantify.

\subsection{Examples and non-examples}

We begin with some examples of definable sets.

\begin{example}
Any sort is a definable set relative to any theory as is any finite product of sorts.  Perhaps less trivially, the range of any term is definable, since if we are given a term $\tau(\bar x)$ and we wish to quantify over the range of $\tau$, for instance $\inf\{ \varphi^\M(y)\colon y \in \text{range}(\tau^\M)\}$, we would evaluate $\inf_{\bar x} \varphi(\tau(\bar x))$.  This fact already pays dividends in the case of operator algebras, where we can express any self-adjoint element as $(x + x^*)/2$ and any positive element as $x^*x$, hence both of these sets of elements are definable relative to theory of \cstar-algebras.
\end{example}

\begin{example}
A slightly less trivial example is that the set of projections is definable relative to the class of \cstar-algebras.  Recall that the element $p$ in a \cstar-algebra $A$ is a projection if it is self-adjoint and satisfies $p^2 = p$.  To see that the set of projections is definable, we use some functional calculus.  Suppose we have a self-adjoint element $b$ in some \cstar-algebra $A$ satisfying $\|b^2 - b\| < \e$, where $\e < 1/4$.  We know that the subalgebra $B$ of $A$  generated by $b$ and $1_A$ is isomorphic to $C(\sigma(b))$ by an isomorphism sending $b$ to the identity function and $1_A$ to the constant function $1$, where $\sigma(b)$ is the spectrum of $b$ and $C(\sigma(b))$ is the \cstar-algebra of complex-valued continuous functions on $\sigma(b)$.  By applying the isomorphism, in $C(\sigma(b))$, we have $\|x^2 -x\| < \e$ and so $\sigma(b)$ can be partitioned into two disjoint open sets $O_0$ and $O_1$ with $i \in O_i$ for $i = 0,1$.  If we let $\pi$ be the function which is constantly $i$ on $O_i$, then this is a continuous function on $\sigma(b)$.  Additionally, $\pi$ is a projection and $\|\pi - x\| < \e$.  Pulling the isomorphism from $C(\sigma(b))$ back into $B$, we find a projection $p \in B$ such that $\|b - p\| < \e$.  By Theorem \ref{def-thm} (3), the set of projections is definable relative to the theory of \cstar-algebras.
\end{example}

We now give two non-examples to demonstrate that the notion of definable set is quite distinct from the notion of zero set.  The first is a somewhat abstract example.

\begin{example}
    The structure $\M$ has a single sort with underlying set $[0,1]$ equipped with the discrete metric and a single relation symbol $R$ which is interpreted as $R(x) = x$. Let $T = Th(\M)$ and consider the $T$-functor $X$, for $\N \in Mod(T)$, defined by $X(\N)=\{ a \in \N \colon R^\N(a) \leq 1/2 \}$.  There are a variety of ways to see that $X$ is not definable.  We will use the ultraproduct characterization.  Notice that if $\cU$ is a non-principal ultrafilter on $\bbN$, then the sequence $(1/2 + 1/n : n \in \bbN)$ has $R$-value $1/2$ in $\M^\cU$ but the distance from this element to every element in $[0,1/2]^\cU$ is 1.  This shows that $X(\M)^\cU \neq X(M^\cU)$ and so $X$ is not $T$-definable.
  \end{example}

  \begin{example}

An example that appears in \cite{munster} in the context of \cstar-algebras is the set of normal elements.  In a \cstar-algebra, the element $a$ is said to be normal if it commutes with its adjoint, hence we are looking at the zero set of $\|a^*a - aa^*\|$.  The set of normal elements in a \cstar-algebra is not definable. (See  \cite[Proposition 3.2.0.8]{munster} for a proof.)
\end{example}

\subsection{Principal types}

For emphasis, we record the following remark made during Comment \ref{definablecomment}.

\begin{prop}\label{def-prop}
If $X$ is a $T$-definable set, then $X(\M)$ is non-empty for some model $\M$ of $T$ if and only if $X(\M)$ is non-empty for all models $\M$ of $T$.
\end{prop}

% \proof By Theorem \ref{def-thm}, there is a formula which when interpreted in $\M$ for any model of $T$ measures the distance to $X(\M)$.  If $X(\M)$ is empty then the distance to $X(\M)$ does not exist (or is infinite).  So $X(\M)$ is non-empty for all models $\M$ of $T$. \qed 

The following proposition makes a connection between the logic topology and the metric on $S_{\bar x}(T)$.  We say that a (partial) type is \emph{principal}\index{principal type} relative to a theory $T$ if the $T$-functor $X$ which sends $\M$ to the set of realizations of $p$ in $\M$ is a $T$-definable set.  By Proposition \ref{def-prop}, any principal type is realized in all models of $T$; the converse of this statement is the continuous form of the omitting types theorem (proven in Section 10 below), and uses the next result in a crucial way.

\begin{prop}\label{principaltype}
For $p \in S_{\bar x}(T)$, the following are equivalent:
\begin{enumerate}
    \item $p$ is principal relative to $T$.
    \item Every $\e$-ball around $p$ contains a logical open neighborhood of $p$.
    \item Every $\e$-ball around $p$ has nonempty interior in the logic topology.
\end{enumerate}
\end{prop}

\proof 
%The equivalence of (1) and (2) is the equivalence of (1) and (4) from Theorem \ref{def-thm}.

If $p$ is principal, then (2) follows from Theorem \ref{def-thm} (4). 

We now prove that (2) implies (1).  This argument is based on the proof of \cite[Theorem 9.12]{BBHU}.
%Suppose that $\M$ is $|\cL|$-saturated which would mean that $p(\M)$ is non-empty.
Fix $n \in \bbN$ and consider a formula $\varphi_n(x)$ such that $p(\varphi_n) = 0$ and $p \in O^T(\varphi_n,(-\infty,\delta_n)) \subseteq B(p,1/n)$.  It follows that if $\M$ is a model of $T$ and $a \in \M$ satisfies $\varphi_n^\M(a) < \delta_n$, then $d(tp^\M(a),p) \leq 1/n$.  Now consider the $T$-formula $P(x)$ defined by
\[
\sum_n \frac{|\varphi_n(x)|}{B_n2^n},
\]
where $B_n$ is the bound for the formula $\varphi_n$ (see Example \ref{weighted}). Notice that, for all $\e>0$, there is $\delta>0$ so that, for all models $\M$ of $T$ and all $a\in \M$, if $P^\M(a)<\delta$, then $d(tp^\M(a),p) <\e$.  Using Proposition \ref{2.10}, we can find an increasing continuous function $\alpha$ with $\alpha(0)=0$ such that
$d(tp(a),p) \leq \alpha(P^\M(a))$ holds for all models $\M$ of $T$ and all $a \in \M$.  Setting $D(x) = \inf_y(\alpha(P(y))+ d(x,y))$, we claim that $D(x) = d(x,p(\M))$ for all models $\M$ of $T$.

To see this, first assume that $\M$ is $\omega$-saturated, which in particular implies that for all $a \in \M$, we have $d(a,p(\M)) = d(tp^\M(a),p)$.  Fix $a \in \M$; we first show that $D^\M(a) \leq d(a,p(\M))$.  Recalling that $P(y) = 0$ whenever $y$ realizes $p$, we see that
\[
D^\M(a) = (\inf_y(\alpha(P(y)) + d(a,y)))^\M \leq \inf_{y \in p(\M)}d(a,y) = d(a,p(\M)).
\]
Note that this inequality did not use the $\omega$-saturation assumption.  For the other inequality, given $b \in \M$, we have $\alpha(P^\M(b)) \geq d(tp^\M(b),p)=d(b,p(\M))$, hence
\[
D^\M(a) = (\inf_y(\alpha(P(y)) + d(a,y)))^\M \geq \inf_{y\in \M}(d(y,p(\M)) + d(a,y)) \geq d(a,p(\M)).
\]
%In the other direction, assume that $d(a,p(\M) = r$ and $\e > 0$.  Choose $b \in p(\M)$ such that $d(a,b) \leq r+ \e$.  But we have $D^\M(a) \leq d(a,b) \leq d(a,p(\M)) + \e.$ Let $\e$ tend to zero and we have $D^\M(a) \leq d(a,p(\M)$.  So in $\M$, $D(x) = d(x,p(\M))$.  

Now suppose that $\N$ is an arbitrary model of $T$ and $\M$ an $\omega$-saturated elementary extension, hence $D^\M(x) = d(x,p(\M))$ by the previous paragraph. We now show that $D^\N(a) = d(a,p(\N))$ for all $a\in \N$.  As noted above, the inequality $D^\N(a) \leq d(a,p(\N))$ holds without the $\omega$-saturation assumption.  For the other inequality, fix $\e > 0$; we will show $d(a,p(\N)) \leq D^\N(a) + \e$, which suffices by letting $\e$ tend to 0.  Towards this end, first notice that if $a \in \N$ is such that $P^\N(a) = 0$, then $a$ realizes $p$ and so $p(\N) = Z(P^\N)$.  

Now choose $a = a_0\in \N$.  We know that $D^\M(a) = d(a,p(\M))$ which implies that there is $c \in \M$ so that $D^\M(c) = 0$ and $d(a,c) = D^\M(a)$.  So $\M$ satisfies $\inf_y \max\{D(y),|d(a,y) - D(a)|\}$. By elementarity, we can choose $a_1 \in \N$ so that $D^\N(a_1) \leq \frac{\e}{8}$ and $|D^\N(a_0) - d(a_0,a_1)| \leq \frac{\e}{8}$. We continue to create a Cauchy sequence $(a_n)$ in $\N$ such that
\begin{itemize}
    \item $D^\N(a_n) \leq \frac{\e}{2^{n+2}}$, and
    \item $|D^\N(a_n) - d(a_n,a_{n+1})| \leq \frac{\e}{2^{n+2}}$.
\end{itemize}
Given that we have constructed $a_n$, we note that since $D^\N(a_n) \leq \frac{\e}{2^{n+2}}$, in $\M$ we can find $c$ realizing $p$ so that $D^\M(c) = 0$ and $D^\M(a_n) = d(a_n,c)$.  By elementarity of $\N$ in $\M$, we can construct $a_{n+1}$ satisfying the properties listed above.
By the completeness of $\N$, the sequence $(a_n)$ has a limit $b \in \N$ and by the continuity of $D^\N$, we have that $D^\N(b) = 0$, hence $b$ realizes $p$.  We also conclude from our two items that $d(a_n,a_{n+1}) \leq \frac{\e}{2^{n+1}}$ for all $n$.  We thus have
\[
%D^\N(a) \leq d(a_0,a_1) + \e/8 \leq d(a_0,a_2) + d(a_2,a_1) + \e/8 \leq d(a_0,b) + \sum_n d(a_{n+1},a_{n+2})
d(a,p(\N))\leq d(a,b) = \lim_n d(a,a_n) \leq \lim_n (d(a_0,a_1) + \sum_{j=1}^{n-1} d(a_j,a_{j+1})) \leq D^\N(a) + \e.
\]
We have thus finished the proof that (2) implies (1).

 The direction (2) implies (3) is trivial, so it remains to prove (3) implies (2).  Suppose that the nonempty logic open set $O^T(\varphi,(r,s))$ is contained in the $\e$-ball around $p$.  Choose $r'$ and $s'$ so that $r < r' < s' < s$ and such that, setting $\psi(x)=\inf_y \max(r' \dotminus \varphi(y), \varphi(y)\dotminus s',d(x,y) \dotminus \e)$, we have $p(\psi)=0$. It follows that $p$ belongs to the logic open set $O^T(\psi,(-\infty,s''))$ for any $s''>0$; if $r<r'-s''$, $s'+s''<s$ and $s''<\epsilon$, then $O^T(\psi,(-\infty,s''))$ is contained in the $3\e$-ball around $p$.  Since $\e$ is arbitrary, this proves (2). \qed

%
%\textbf{Another argument for (2) implies (1)}:  Fix $n \in \bbN$ and consider a formula $\varphi_n(x)$ such that $p(\varphi_n) = 0$ and $p \in O^T(\varphi_n,(-\infty,\delta_n)) \subseteq B(p,2^{-n})$.  Fix $\M\models T$ and suppose $a\in \M$ is such that $\varphi_n^\M(a)<\delta_n$; we claim that $d(a,p(\M))\leq 2^{-n}$.  Fix $m\geq 1$.  In an elementary extension $\N$ of $\M$, there is a realization $b$ of $p$ with $d(a,b)\leq 2^{-n}$.  Consequently, there is $a_1\in \M$ such that $d(a,a_1)\leq \frac{1}{2^{n}}+\frac{1}{2^{m+1}}$ and $\varphi_{q_1}^\M(a)<\delta_{q_1}$ for some $q_1$ sufficiently large.  By the same reasoning, there is some $q_2>q_1$ sufficiently large and $a_2\in \M$ such that $d(a_1,a_2)\leq \frac{1}{2^{m+2}}$ and $\varphi_{q_2}^\M(a_2)<\delta_{q_2}$.  In this way, we obtain a Cauchy sequence $(a_k)$ in $\M$ and a sequence $\delta_{q_k}$ such that $\lim_k \delta_{q_k}=\infty$, $d(a_k,a_{k+1})\leq \frac{1}{2^{m+k+1}}$, and with $\varphi^\M_{q_k}(a_k)<\delta_{q_k}$.  Let $a'\in \M$ denote the limit of $(a_k)$.  It follows that $d(a,a')\leq \frac{1}{2^n}+\frac{1}{2^m}$.  Moreover, letting $\epsilon_k:=d(a,a_k)$, we then have that there is some $b_k$ in an elementary extension of $\M$ realizing $p$ such that $d(a,b_k)\leq \epsilon_k+2^{-q_k}$.  It follows that $a'$ realizes $p$.  Since $m$ is arbitrary, we see that $d(a,p(\M)\leq 2^{-n}$, as desired.

\section{Quantifier elimination}\label{sec-qe}

We say that a theory $T$ in a language $\cL$ has \emph{quantifier elimination}\index{quantifier elimination} if, for every $\cL$-formula $\varphi(\bar x)$, there is a quantifier-free $\cL$-formula $\psi(\bar x)$ which is $T$-equivalent to $\varphi(\bar x)$.  In terms of basic $\cL$-formulas, this says: for every $\e > 0$ and every basic $\cL$-formula $\varphi(\bar x)$, there is a quantifier-free basic $\cL$-formula $\psi$ such that $\|\varphi(\bar x) - \psi(\bar x)\|_T < \e$.

We will give three characterizations of quantifier elimination.  The first is syntactical.

\begin{prop}\label{syn-qe}
Suppose that $T$ is an $\cL$-theory. Then $T$ has quantifier-elimination if and only if every formula of the form $\inf_x \varphi(x,\bar y)$, where $\varphi(x,\bar y)$ is a basic $\cL$-formula, is $T$-equivalent to a quantifier-free $\cL$-formula.
\end{prop}

\proof 
The forward direction is immediate and the backward direction is proven by induction on the formation of basic $\cL$-formulas, the quantifier case being taken care of by the backwards assumption and the connective case being immediate.
% We only need to prove right implies left. From the assumption, it is immediate, by induction on the formation of basic formulas, that every basic formula is equivalent to a quantifier-free formula.  The only issue is when one is considering a formula of the form $\inf_x \varphi(x,\bar y)$ where $\varphi(x,\bar y)$ is a limit of quantifier-free basic formulas.  However, since the approximation of $\varphi(x,\bar y)$ by quantifier-free formulas is uniform, $\inf$ and the limit are interchangeable so if $\varphi = \lim_{n \rightarrow \infty} \varphi_n$ where $\varphi_n$ are basic quantifier-free formulas then
% $\inf_x \varphi(x,\bar y)$ is $T$-equivalent to $\lim_{n \rightarrow \infty} \inf_x \varphi_n(x,\bar y)$ and so by the operative assumption, the latter is a quantifier-free formula. 
\qed

\begin{prop}\label{qftype-qe}
Suppose that $T$ is an $\cL$-theory. Then $T$ has quantifier-elimination if and only if whenever $\M$ is a model of $T$ and $\bar a,\bar b \in \M$ have the same values for basic atomic formulas, then they satisfy the same complete type. 
\end{prop}

\proof This follows from Proposition \ref{gen-SW} applied to the real algebra system of quantifier-free formulas. \qed

\begin{example}
We can use Proposition \ref{qftype-qe} to see that the theory $T$ of the Urysohn sphere $U$ has quantifier elimination.  Indeed, suppose that $\N$ is a model of $T$ and $\bar a,\bar b \in \N$ have the same quantifier-free type.  Without loss of generality, we can assume that $\N$ is separable and thus we may further assume that $\N$ is $U$ itself by the back-and-forth argument given in Example \ref{ury-example}.  Since $\bar a$ and $\bar b$ have the same quantifier-free type, the map $f$ sending $\bar a \mapsto \bar b$ is an isometry.  Since $U$ is ultrahomogeneous, $f$ extends to an automorphism of $U$, which implies that $\bar a$ and $\bar b$ have the same complete type in $U$, as desired.
\end{example}

A more semantic characterization of quantifier elimination is given by the following:

\begin{prop}
Suppose $T$ is an $\cL$-theory. Then $T$ has quantifier elimination if and only if $T$ has the property that whenever $\M$ is a model of $T$, $\bar a$ and $\bar b \in \M$ have the same quantifier-free type, then there is $\N$ with $\M \prec \N$ such that the map sending $\bar a$ to $\bar b$ can be extended to an embedding of $\M$ into $\N$. 
\end{prop}

Note that if $\M$ in the Proposition is separable then $\N$ could be taken to be $\M^\cU$ where $\cU$ is a non-principal ultrafilter on $\bbN$.

\proof We leave the proof of the forward direction as an exercise to the reader and sketch a proof of the backwards direction using Proposition \ref{syn-qe}.  Fix a basic quantifier-free $\cL$-formula $\varphi(x,\bar y)$. Suppose that for some $\e$, the following set of sentences is satisfiable:
\begin{multline*}
\{ |\psi(\bar c) - \psi(\bar d)| \colon \psi(\bar y) \text{ is a basic quantifier-free $\cL$-formula}\}  \\
 \cup \{\e \dotminus |\inf_x \varphi(x,\bar c) - \inf_x \varphi(x,\bar d)|  \},
\end{multline*}
where $\bar c$ and $\bar d$ are sequences of new constants.  It follows that there is a model $\M$ of $T$ and $\bar a,\bar b\in \M$ with the same quantifier-free type but $(\inf_x \varphi)^\M(x,\bar a) < (\inf_x \varphi)^\M(x,\bar b)$.  By assumption, there is an elementary extension $\N$ of $\M$ such that the map $\bar a\mapsto \bar b$ extends to an embedding $f$ from $\M$ to $\N$.  However, if we choose $c \in\M$ such that $\varphi^\M(c,\bar a) < (\inf_x\varphi)^\M(x,\bar b)=(\inf_x \varphi)^\N(x,\bar b)$, then we have a contradiction when we consider the value of $\varphi^\N(f(c),\bar b)$.

Consequently, the above set of sentences cannot be satisfiable for any $\e$.  By compactness, for any $\e>0$, we can find, after some manipulation, a quantifier-free basic formula $\psi(\bar x)$ and $\delta>0$ such that, for any model $\M$ of $T$ and $\bar a, \bar b \in \M$,
\[
|\psi(\bar a) - \psi(\bar b)| < \delta \text{ implies } |\inf_x \varphi(x,\bar a) - \inf_x \varphi(x,\bar b) | \leq\e.
\]
By Proposition \ref{gen-SW} applied to the real algebra system of quantifier-free formulas, we have that $\inf_x \varphi(x,\bar y)$ is equivalent to a quantifier-free formula.  By Proposition \ref{syn-qe}, $T$ has quantifier elimination.
\qed
\begin{exercise}
Prove the forward direction of the previous proposition.
\end{exercise}

\begin{example}
We can use the previous proposition to show that theory of in-finite-dimensional Hilbert spaces admits quantifier-elimination.  (We leave it to the reader as an exercise to verify that infinite-dimensionality is expressible in continuous logic.)  Indeed, suppose that $\cal H$ is an infinite-dimensional Hilbert space and that $\bar a, \bar b\in \cal H$ are tuples with the same quantifier-free type.  Letting $\cal H_0$ and $\cal H_1$ denote the finite-dimensional subspaces of $\cal H$ generated by $\bar a$ and $\bar b$ respectively. It follows that the map $a_i\mapsto b_i$ for $i=1,\ldots,n$ extends to an isomorphism $\cal H_0\to \cal H_1$.  Since $\cal H$ is infinite-dimensional, we have that the dimension of the orthogonal complement $\cal H_0^\perp$ of $\cal H_0$ in $\cal H$ is the same as the dimension of the orthogonal complement $\cal H_1^\perp$ of $\cal H_1$ in $\cal H$.  Thus, mapping $\cal H_0^\perp$ isomorphically into $\cal H_1^\perp$ allows one to extend the original mapping to an automorphism of $\cal H$.  Note that in this example, there is no need for an elementary extension; in more complicated examples, passing to an elementary extension is in fact necessary.
\end{example}

\section{Imaginaries}\label{sec-imag}
In classical model theory, imaginary elements play an important role.  The construction of the theory $T^{eq}$ from $T$ creates the largest expansion of $T$ for which the category of models is equivalent.  Although somewhat more complicated in the continuous setting, one can define imaginaries for metric structures which have similar properties.  We consider them in the next three subsections.

\subsection{Products}
Fix a continuous language $\cL$ and an $\cL$-theory $T$.  Suppose that $\bar S = (S_n : n \in \bbN)$ is a sequence of sorts from $\cL$.  We add one new sort $S$ to $\cL$ to form $\cL'$.  The intention is that we will expand models of $T$ to include as the interpretation of $S$ the product of the sorts $S_n$.  We expand every model $\M$ of $T$ by interpreting $S$ as $\prod_{n \in \bbN} S_n^\M$.  We need a metric for this sort.  The choice of metric is not unique but for definiteness we use, for $\bar x, \bar y \in \prod_{n \in \bbN} S_n^\M$,
\[
d(\bar x,\bar y) = \sum_{n \in \bbN} \frac{d_n^\M(x_n,y_n)}{B_n2^n},
\]
where $d_n$ is the metric symbol for $S_n$ and $B_n$ is the bound on that metric.  The key property that $d$ has is that for every $\e > 0$, there is an $n\in \bbN$ such that if $\bar x, \bar y \in \prod_{n \in \bbN} S_n^\M$ and $x_i = y_i$ for every $i < n$, then $d(\bar x,\bar y) < \e$.  The only other symbols we add to the language $\cL'$ are maps $\pi_n$ from $S$ to $S_n$, and in the canonical expansion of an $\cL$-structure to an $\cL'$-structure, we interpret $\pi_n$ as the projection onto the $n^{th}$ coordinate.  The uniform continuity modulus for $\pi_n$ can be taken to be $\e/B_n2^n$.

We record here three families of sentences that hold in all expanded structures.
First, for any $n \in \bbN$, we want to state that 
\[
\sum_{i \leq n} \frac{d_i^\M(x_i,y_i)}{B_i2^i} \leq d(x,y) \leq \sum_{i \leq n} \frac{d_i^\M(x_i,y_i)}{B_i2^i} + \frac{1}{2^n}
\]
holds for all $x,y\in S$; we can express this fact in continuous logic by requiring that the following two sets of sentences evaluates to $0$:
\[
\sup_{x,y \in S} \left ( d(x,y) \dotminus (\sum_{i \leq n} \frac{d_i^\M(\pi_i(x),\pi_i(y))}{B_i2^i} + \frac{1}{2^n}) \right ).
\] and
\[
\sup_{x,y \in S} \left ( \sum_{i \leq n} \frac{d_i^\M(\pi_i(x),\pi_i(y))}{B_i2^i} \dotminus d(x,y) \right ).
\]
We would also like the finite products to be dense in $S$.  We can achieve this by satisfying the following sentences, one for each $n\in \bbN$:
\[
\sup_{x_i \in S_i, i < n}\inf_{x \in S} \max_{i < n} \{d_i(\pi_i(x),x_i)\}.
\]

Let $T'$ be the $\cL'$-theory consisting of $T$ together with the three preceding sets of sentences.  If $\M'$ is a model of $T'$, then it has the form $(\M,Y)$, where $\M$ is a model of $T$ and $Y$ is the interpretation of the sort $S$.  We wish to see that $S$ is essentially $\prod_{n \in \bbN} S_n^\M$.  Fix $x_i \in S_i^\M$.  The third set of axioms guarantees that, for every $n \in \bbN$, there is $y^n \in Y$ such that $d(\pi_i(y^n),x_i) < \frac{1}{2^n}$ for every $i < n$.  The first set of axioms assures us that the sequence $(y^n : n \in \bbN)$ converges to some $y \in Y$ (since $Y$ is complete).  For this $y$, we have that $d(\pi_n(y),x_n) = 0$ for all $n\in \bbN$.  If there is another $y' \in Y$ with the same property, then using the first set of sentences again, we have that $d(y,y') = 0$ and so $y = y'$.  In this way we see that $Y$ can be identified with $\prod_{n \in \bbN} S_n^\M$ via the projection maps.  This identification is an isometry since the first two sets of axioms guarantees that the metric on $Y$ is given by
\[
d(x,y) = \sum_n \frac{d_n(\pi_n(x),\pi_n(y))}{B_n2^n}.
\]
What we have achieved is the following:

\begin{thm}
The forgetful functor from $Mod(T')$ to $Mod(T)$ is an equivalence of categories.
\end{thm}

\subsection{Definable Sets}

Although the following collection of imaginary elements is not entirely essential, it is often useful in practice.  Suppose that $X$ is a definable $T$-functor for some $\cL$-theory $T$.  We create a new language $\cL'$ by adding a new sort $S$ and a new function symbol $i$ with domain $S$ and co-domain equal to the domain of $X$.  The intention is to interpret $S$ as $X(\M)$ for $\M$ a model of $T$ and to interpret $i$ as the inclusion map from $S$ into $\M$.  Here are sentences in $\cL'$ which achieve what we want:
\[
\sup_{x \in S} d(i(x),X)
\]
and
\[
\sup_{x,y \in S} |d(x,y) - d(i(x),i(y))|.
\]
It is reasonably easy to see that the $\cL'$-theory $T'$ formed from $T$ and these two axioms yields an essentially unique expansion of models of $T$.  Once again we have: 
\begin{thm}
The forgetful functor from $Mod(T')$ to $Mod(T)$ is an equivalence of categories.
\end{thm}

\subsection{Canonical Parameters}

The following collection of imaginary elements is the closest analogy to the classical case.  Fix a language $\cL$ and an $\cL$-theory $T$. Choose an $\cL$-formula $\varphi(\bar x,\bar y)$.  Suppose that the sort of $\bar y$ is $\bar S$.  The idea is that we wish to capture in a single sort all of the functions of the form $\varphi^\M(\bar x,\bar a)$ for $\M$ a model of $T$ and  $\bar a \in \M$. We will call the elements of this sort the \emph{canonical parameters}\index{canonical parameters}
\footnote{The name comes from the fact that $\bar a$ are the parameters for the function $\varphi(\bar x,\bar a)$ but there may be other parameters which give the same function.  If we make any two such sequences equivalent, then we have a "canonical" choice for the parameters of this function.}
of $\varphi(\bar x,\bar y)$.  In order to do this, we create a new language $\cL'$ by adding to $\cL$ a sort $S$ with metric symbol $d_S$ and a function symbol $\pi$ with domain $\bar S$ and co-domain $S$. We add the following two axioms to the theory $T$ to form $T'$:

\[
\sup_{z \in S}\inf_{\bar y \in \bar S} d(z,\pi(\bar y))
\]
and
\[
\sup_{\bar y,\bar y' \in \bar S} \left | d_S(\pi(\bar y),\pi(\bar y')) - \sup_{\bar x} | \varphi(\bar x,\bar y) - \varphi(\bar x,\bar y')| \right |.
\]

The first axiom says that the range of $\pi$ is dense in $S$ and the second explicitly defines the metric on the range of $\pi$.
If $\M$ is a model of $T$, then it can be expanded to a model of $T'$ by letting $S$ be the completion of $\bar S^\M$ quotiented by $d'$ where
%the set of functions of the form $\varphi^\M(\bar x, \bar a)$ for all $\bar a \in \bar S^\M$ quotiented by the pseudo-metric $d_S$ on $S$ determined by
\[
d'(\bar a,\bar b) = (\sup_{\bar x} | \varphi(\bar x,\bar a) - \varphi(\bar x,\bar b)|)^\M.
\]
We define $\pi(\bar a) = \bar a/d'$ for $\bar a \in \bar S^\M$.  Once again there is essentially only one way to expand a model of $T$ to a model of $T'$ and so we leave it as an exercise to show that:
\begin{thm}
The forgetful functor from $Mod(T')$ to $Mod(T)$ is an equivalence of categories.
\end{thm}

\subsection{Conceptual completeness}
If one iterates the constructions from the previous three subsections to a given $\cL$-theory $T$, one can construct a theory $T^{eq}$ (referred to as $T^{meq}$ in the article by Berenstein and Henson in this volume)  in an expanded language $\cL^{eq}$ which is closed under taking finite and countable products of sorts, and adding definable sets and canonical parameters. (In fact, one can be economical and first close under products, then add definable sets and finally add canonical parameters.)  

In order to state the main theorem regarding $T^{eq}$, we need to say what interpretation means in continuous logic.

\begin{defn}
Let $\cL$ and $\cL'$ be continuous languages and suppose that $T$ and $T'$ are theories in $\cL$ and $\cL'$ respectively. An \emph{interpretation}\index{interpretation}
$I$ of $T$ in $T'$ assigns:
\begin{enumerate}
\item to every sort $S$ of $\cL$ , a sort $I(S)$ of $\cL'$,
\item to every function symbol $f\colon \prod S_i \rightarrow S$ an $\cL'$-formula $I(f)$ with free variables in $\prod S_i \times S$, such that 
\[
T' \models \text{``}I(f)(\bar x, y) = 0 \text{ is the graph of a total function of } \bar x\text{''}
\]
\item to every relation symbol $R$ with domain $S_1\times \cdots \times S_n$ in $\cL$, an $\cL'$-formula $I(R)$ with domain $I(S_1) \times \cdots \times I(S_n)$,
\end{enumerate}
such that, for any $\cL$-sentence $\varphi$, $T \models \varphi$ iff $T' \models I(\varphi)$ (where by $I(\varphi)$ one means the result of applying the interpretation of the language inductively to $\varphi$).
\end{defn} 

While this definition is not the most general definition of interpretation, it is sufficient to prove the following theorem, whose proof can be found in \cite{AH}.

\begin{thm}[Conceptual completeness]
Suppose that $T$ is an $\cL$-theory.  Then $T^{eq}$ satisfies the following properties:
\begin{enumerate}
\item The forgetful functor from $Mod(T^{eq})$ to $Mod(T)$ is an equivalence of categories.
\item If $T'$ is an expansion of $T$ in some language $\cL'$ such that the forgetful functor from $Mod(T')$ to $Mod(T)$ is an equivalence of categories, then $T'$ can be interpreted in $T^{eq}$.
\end{enumerate}
\end{thm}

\section{Omitting types}\label{sec-omit}

We say that a structure $\M$ \emph{omits}\index{omitting type} a (partial) type $p$ if $p$ is not realized in $\M$. Criteria in the classical case for the omission of types is well understood.  The situation in the continuous case is more problematic; see the discussion at the end of this section.  We will see here that the role of definability is pivotal in the omitting types theorem for continuous logic.  

Throughout this section, we work in a theory $T$ in a separable language $\cL$.  As mentioned before, by Proposition \ref{def-prop}, if $p$ is omitted in some model, then $p$ is not principal. The converse is the omitting types theorem for continuous logic:
\begin{thm}[Omitting types theorem]
For each $n\in \bb N$, suppose that $p_n$ is a complete type which is not principal.  Then there is a separable model of $T$ which omits each $p_n$.
\end{thm}

\proof This proof provides an opportunity to show what a Henkin construction looks like in the continuous setting. We will sketch the proof where we omit a single complete type $p$ in one variable and leave it to the reader to handle the more general case. To this end, we consider the language $\cL'$ which is $\cL$ together with countably many new constants symbols $c_n$ for $n \in \bbN$.  The intention is that we will construct an $\cL'$-theory $T'$ extending $T$ and a model of $T'$ in which the interpretations of the $c_n$'s are dense.  
Recall the notation for open sets in the logic topology: a basic open set relative to the theory $T$ has the form
\[
 \{ \bar T \colon T \subseteq \bar T \text{ is a complete $\cL'$-theory and } \bar T \models a < \psi < b \}
\]
for some $\cL'$-sentence $\psi$ and interval $(a,b)$.  We denote this open set by $O^T(\psi,(a,b))$.

We construct $T'$ inductively as the intersection of a decreasing sequence of non-empty basic logical open sets $O_n$ for $n \in \bbN$.  We set $O_0$ to be the set of all complete $\cL'$-theories extending $T$.  Now, fix a countable dense set of $\cL'$-sentences $\psi_n$ for $n \in \bbN$ and a countable dense set of $\cL'$-formulas in one free variable $\varphi_n(x)$ for $n \in \bbN$.  We will want the set $O_n$ to satisfy the following conditions:

First, for every $m$ and every rational number $r$, there is $n\in \bb N$ and rational interval $(a,b)$ with $b-a < r$ such that 
\[
O_n \subseteq  O^T(\psi_m,(a,b)).
\]

Second, for every $m\in \bb N$ and rational number $r$, if for some $n$ we have
\[
O_n \subseteq O^T(\inf_x \varphi_m(x),(-\infty,r)),
\]
then for some $k > n$ and $l\in \bb N$ we have
\[
O_k \subseteq O^T(\varphi(c_l),(-\infty,r)).
\]
It is not hard to see that if we can succeed in obtaining these properties, then $\bigcap_{n \in \bbN} O_n$ is a complete $\cL'$-theory such that in any model, the closure of the interpretations of the constants $c_n$ is the underlying set of a model of $T$.

In order for the aforementioned model of $T$ to omit $p$, we need to demand more from the open sets $O_n$.  In the classical case, it is enough to require that none of the new constants satisfy $p$.  Here, since the new constants are merely dense in the model of $T$ we have constructed, we will instead need to keep the constants some fixed distance away from the set of realizations of $p$.  To this end, recall that since $p$ is not principal, by Proposition \ref{principaltype}, there is some $\e > 0$ such that the open $\e$-ball around $p$ does not contain any logical open set.  Suppose that we are at some stage of the construction and we have constructed the logical open set $O_n$, which we may assume, after some manipulation, has the form $O^T(\theta(c_k,\bar d), (-\infty,0))$, where $c_k$ is the new constant we are currently interested in keeping away from $p$ and $\bar d$ is the rest of the new constants in $\theta$.  
Now consider the logical open set $O^T(\inf_{\bar y}\theta(x,\bar y), (-\infty,0))$.  By the assumption that $p$ is not principal, there is some $q \in S_x(T)$ such that $q$ belongs to $O^T(\inf_{\bar y}\theta(x,\bar y),(-\infty,0))$ and $d(p,q) > \frac{\e}{2}$.  By compactness, we can find some $\cL$-formula $\psi(x)$ such that:
\begin{itemize}
   % \item $q$ belongs to the logical open set $O^T(\psi,(-\infty, 0))$,
    \item for all $\M\models T$ and $a\in \M$ with $\psi^\M(a)<0$, we have $(\inf_{\bar y}\theta(a,\bar y))^\M<0$, and
    \item for any $\M\models T$ and  $a,b\in \M$ such that $a$ realizes $p$ and $\psi^\M(b)<0$, we have $d(a,b) > \frac{\e}{2}$. 
    \end{itemize}
We now let $O_{n+1} = O_n \cap O^T(\psi(c_k),(-\infty,0))$ and note that this will guarantee that no interpretation of $c_k$ will be within $\frac{\e}{2}$ of a realization of $p$ in any model of $T$.  Since we can achieve this for all of the new constants, no Cauchy sequence of the new constants will approach a realization of $p$ and so the model we construct will omit $p$. \qed
%  Let $\Sigma_0 = T$.  We will want to create an increasing chain of $\cL'$-theories $\Sigma_n$ such that
%\begin{enumerate}
%\item For each $n$, $\Sigma_n \setminus T$ is finite.
%\item For every $\psi_m$ and rational interval $(a,b)$ there is $n$ such that either $\psi_m \in (a,b)$ is in $\Sigma_n$ or $\psi_m \not \in (a,b)$ is in $\Sigma_n$.
%\item For every $\varphi_m(x)$ and rational $r$, if for some $n$ $\inf_x \varphi(x) < r$ is in $\Sigma_n$ then for some $k >n$ and some $c_l$,
%$\varphi(c_l) < r$ is in $\Sigma_k$.
%\end{enumerate}
%These three conditions are the continuous analogs of the usual conditions involved in a Henkin construction.  In order to guarantee that the model we construct omits $p$, it will not be enough to guarantee that each constant does not satisfy $p$.  We will need to know that the realization of each constant stays some fixed distance away from realizations of $p$.  Since $p$ is not definable, there is some $\e > 0$ such that in $S_x(T)$, the $\e$-ball around $p$ does not contain any logical open set. So suppose that $\Sigma_n = T \cup \varphi(c,\bar d)$

The situation regarding omitting partial types in continuous logic is quite problematic as evidenced by the results of Farah and Magidor from \cite{FM}.  There they show the following:
\begin{enumerate}

\item There is a complete theory $T$ in a separable language and countably many partial types such that every finite subset is omitted in some model of $T$ but no model of $T$ simultaneously omits all of them.
\item There is a complete theory $T$ in a separable language and partial types $s$ and $t$, each one omitted in some model of $T$, but no model omits both.
\end{enumerate}

A variant of the omitting types theorem used to omit certain kinds of partial types has proven useful in the operator algebra context; for example, see \cite{munster}, \cite{enforceable}, and Goldbring's article in this volume.

\section{Separable categoricity}\label{sec-sep-cat}

For an infinite cardinal $\lambda$, we say that a theory $T$ is \emph{$\lambda$-categorical}\index{categorical} if $\lambda \geq \chi(T,\cL)$ and whenever $\M$ and $\N$ are models of $T$ such that $\chi(\M) = \chi(\N) = \lambda$, then $\M \cong \N$.  We will focus on the case where $\lambda$ is $\aleph_0$ and we will typically say \emph{separably categorical}\index{separably categorical} instead of $\aleph_0$-categorical.  As in the classical case, we have the following theorem.

\begin{thm}\label{RN}
Suppose $T$ is an $\cL$-theory which has only models of density character at least $\chi(T,\cL)$ and is $\lambda$-categorical.  Then $T$ is complete.
\end{thm}

\proof Suppose that $\M$ and $\N$ are models of $T$.  If $\M$ has density character at least  $\lambda$, then by the Downward L\" owenheim-Skolem theorem, we can find $\M_0 \prec \M$ with $\chi(\M_0) = \lambda$.  If $\chi(\M) < \lambda$, then by the upward L\" owenheim-Skolem theorem applied to $Elem(\M)$, we can find $\M_0$ with $\chi(\M_0) = \lambda$ and $\M \prec \M_0$.  Applying the same reasoning to $\N$, we can find $\M_0$ and $\N_0$ so that
\[
\M \equiv \M_0 \cong \N_0 \equiv \N.
\]
From this we conclude that $T$ is complete. \qed

We now state the Ryll-Nardewski Theorem for continuous logic; the proof, done originally in the positive bounded case, is due to Henson.

\begin{thm}
Suppose that $\cL$ is separable and $T$ is an $\cL$-theory with no compact models.  Then $T$ is separably categorical if and only $T$ is complete and every complete type is principal.
\end{thm}

\proof If $T$ is separably categorical then it is complete.  If there is a complete type $p$ which is not principal, then by the omitting types theorem, there is some separable model of $T$ which omits $p$.  However, there is also a separable model of $T$ which realizes $p$; since these two models of $T$ cannot be isomorphic, $T$ is not separably categorical; contradiction.

Now suppose that $T$ is complete and every complete type is principal.  Consider two separable models $\M$ and $\N$ of $T$.  We need to show that $\M \cong \N$.  

We construct two sequences
$a_0^0, a^1_0a^1_1, a^2_0a^2_1a^2_2, \ldots$ in $\M$ and
$b_0^0, b^1_0b^1_1, b^2_0b^2_1b^2_2, \ldots$ in $\N$ such that
all initials segments of the same length have the same type, that is, for any $k$, there is a fixed type for
$a^n_0\ldots a^n_k$ and $b^n_0 \ldots b^n_k$ independent of $n \geq k$.
We will also want to arrange that for every $k$, $\langle a^n_k : n \geq k \rangle$ and $\langle b^n_k : n \geq k \rangle$ form Cauchy sequences converging to $\bf a_k$ and $\bf b_k$ respectively so that 
$\{ \bf a_k : k \in \bbN \}$ and $\{ \bf b_k : k \in \bbN \}$ are dense in $\M$ and $\N$ respectively.
If we can achieve this, then the map sending $\bf a_k$ to $\bf b_k$ extends to an isomorphism from $\M$ to $\N$.

To start, we enumerate countable dense subsets in $\M$ and $\N$; call them respectively $\langle c_k : k \in N\rangle$ and $\langle d_k : k \in N\rangle$.
At stage 0, let $a^0_0 = c_0$.  By assumption, $tp^\M(c_0)$ is principal and hence realized in $\N$ by some $b^0_0$. In general, at each step we alternate by either choosing a $c_k$ or $d_k$ and we revisit each $c_k$ and $d_k$ infinitely often in the construction.
Assume we have chosen $a^n_0\ldots a^n_n$ already and we consider whatever $c_k$ is given to us at this stage.
Let $p(x_0,\ldots,x_n)$ be the type of $a^n_0\ldots a^n_n$ in $\M$ and $q(x_0,\ldots,x_{n+1})$ be the type of $a^n_0\ldots a^n_nc_k$ in $\M$.
Since $q$ is principal, we have a formula $d_q(x_0,\ldots,x_{n+1})$ expressing the distance to realizations of $q$.
Since $d_q^\M(a^n_0\ldots a^n_n,c_k) = 0$, we have that  $\inf_y d_q^M(a^n_0\ldots a^n_n,y) = 0$.
Since $b^n_0\ldots b^n_n$ satisfies $p$ by assumption, we have $\inf_y d_q^\N(b^n_0\ldots b^n_n,y) = 0$.
This means we can find $b^{n+1}_0\ldots b^{n+1}_{n+1}$ realizing $q$ and such that $d(b^n_i,b^{n+1}_i) \leq 1/2^n$ for $i = 0,\ldots,n$.
This guarantees we have the required Cauchy sequences and we have the required density as well. To see the latter fact, fix $\epsilon > 0$. Choose $N$ large enough so that $\sum_{n \geq N} 1/2^n < \epsilon$.  If we visit $c_k$ at stage $t > N$,
then $\bf a_t$ is within $\epsilon$ of $c_k$ and so the $\bf a_k$'s are dense in $M$.  Similarly, the $\bf b_k$'s are dense in $N$, finishing the proof.
 \qed

 Note that the condition ``every complete type is principal'' in Theorem \ref{RN} is equivalent to the statement that the logic and metric topologies on $S_{\bar x}(T)$ agree for all finite tuples $\bar x$ of variables.

 \begin{exercise}
Prove that a complete theory $T$ in a separable language with no compact models is separably categorical if and only if $S_{\bar x}(T)$ is compact in the metric topology for all finite tuples $\bar x$ of variables.
 \end{exercise}

 There are a number of nontrivial examples of separably categorical continuous theories:

\begin{example}
The following theories are separably categorical:
\begin{enumerate}
    \item The theory of infinite-dimensional Hilbert spaces.
    \item The theory of the Urysohn sphere.
    \item The theory of atomless probability algebras.
\end{enumerate}
\end{example}

The fact that the theory of infinite-dimensional Hilbert spaces is separably categorical follows from the fact that any two separable infinite-dimensional Hilbert spaces are isomorphic to $\ell^2$ while the separable categoricity of the theory of the Urysohn sphere follows from the argument presented in Example \ref{ury-example}.  The separable categoricity of the theory of atomless probability algebras is discussed in Berenstein and Henson's article in this volume \cite{BH}.

In the context of operator algebras, there are very few examples of separably categorical theories.  For instance, for any separable II$_1$ factor $N$, there are continuum many nonisomorphic separable II$_1$ factors elementarily equivalent to $N$ (see Theorem 4.1 in the author's article with Goldbring in this volume).  One notable example of a separably categorical \cstar-algebra is the algebra $C(X)$ for $X$ a totally disconnected compact Hausdorff space; in this case, $C(2^\bb N)$ is the unique separable model.  The proof of this fact follows from Stone duality together with the fact that there is a unique countable atomless Boolean algebra.

We end this section with some related considerations.  We say that a model $\M$ of $T$ is \emph{atomic}\index{atomic model} if $tp^\M(\bar a)$ is principal for every finite tuple $\bar a\in \M$.  Theorem \ref{RN} (and its proof) yields the following:

\begin{thm}

\

\begin{enumerate}
\item Any two separable atomic models of $T$ are isomorphic.
\item If $T$ is a complete theory in a separable language, then $T$ is separably categorical if and only if all separable models of $T$ are atomic.
\end{enumerate}
\end{thm}

In the remainder of this section, for simplicity, we assume that $T$ is a complete theory in a separable language.

We say that a model of $T$ is \emph{prime}\index{prime model} if it embeds elementarily into all other models of $T$.

\begin{exercise}
Prove that a model $\M$ of $T$ is prime if and only if $\M$ is separable and atomic.  Moreover, $T$ can have at most one prime model up to isomorphism.  (Hint:  To show that a separable atomic model of $T$ is prime, use a ``forth only'' version of the back and forth argument appearing in the proof of Theorem \ref{RN}.)
\end{exercise}

\begin{exercise}
Prove that $T$ has a prime model if and only if the principal types are dense in $S_{\bar x}(T)$ for all finite tuples $\bar x$ of variables.
\end{exercise}

We say that $T$ is \emph{small}\index{small theory} if $S_{\bar x}(T)$ is separable in the metric topology for all finite tuples $\bar x$ of variables. The following is Proposition 1.17 in \cite{BYU-dfin}:

\begin{prop}
If $T$ is small, then $T$ has a prime model.
\end{prop}

\printindex
\end{document}